\documentclass[journal,twoside,web]{ieeecolor}

\usepackage[T1]{fontenc}
\usepackage{generic}
\usepackage{cite}
\usepackage{amsmath,amssymb,amsfonts}
\usepackage{algorithmic}
\usepackage{graphicx}
\usepackage{textcomp}
\def\BibTeX{{\rm B\kern-.05em{\sc i\kern-.025em b}\kern-.08em
		T\kern-.1667em\lower.7ex\hbox{E}\kern-.125emX}}
\usepackage{xcolor}

\newtheorem{thm}{Theorem}
\newtheorem{rem}[thm]{Remark}
\newtheorem{defn}[thm]{Definition}

\newtheorem{lem}[thm]{Lemma}

\newtheorem{exa}[thm]{Example}

\begin{document}
\title{Temporal Parallelisation of the HJB Equation and Continuous-Time Linear Quadratic Control}
\author{Simo S\"arkk\"a,~\IEEEmembership{Senior Member,~IEEE,}
        \'Angel F. Garc\'ia-Fern\'andez
\thanks{S. S\"arkk\"a is with the Department
        of Electrical Engineering and Automation, Aalto University, 02150 Espoo, Finland (email: simo.sarkka@aalto.fi).}
\thanks{A. F. Garc\'ia-Fern\'andez is with the Department of Electrical Engineering and Electronics, University of Liverpool, Liverpool L69 3GJ, United Kingdom (email: angel.garcia-fernandez@liverpool.ac.uk). }}

\maketitle

\begin{abstract}
This paper presents a mathematical formulation to perform temporal parallelisation of continuous-time optimal control problems, which can be solved via the Hamilton--Jacobi--Bellman (HJB) equation. We divide the time interval of the control problem into sub-intervals, and define a control problem in each sub-interval, conditioned on the start and end states, leading to conditional value functions for the sub-intervals. By defining an associative operator as the minimisation of the sum of conditional value functions, we obtain the elements and associative operators for a parallel associative scan operation. This allows for solving the optimal control problem on the whole time interval in parallel in logarithmic time complexity in the number of sub-intervals. We derive the HJB-type of backward and forward equations for the conditional value functions and solve them in closed form for linear quadratic problems. We also discuss numerical methods for computing the conditional value functions. The computational advantages of the proposed parallel methods are demonstrated via simulations run on a multi-core central processing unit and a graphics processing unit. 
\end{abstract}


\section{Introduction}
\label{sec:introduction}

Continuous-time optimal control refers to methods that aim at controlling a state-trajectory $x(t)$ of a continuous-time dynamic system using a control signal $u(t)$ such that a given cost functional $C[u]$ is minimised \cite{Lewis_book12,Kirk:2004}. The term ``continuous-time'' in this context refers to time being defined on some continuous interval $[t_0,t_f] \subset \mathbb{R}$ instead of at discrete time instants. Examples of such optimal control problems include the optimal steering of a robot, car, airplane, or spacecraft to follow a given trajectory while minimizing, for example, energy, fuel consumption, or time \cite{LaValle_book06, Athans63}. Other examples include controlling of the behaviour of a chemical process \cite{Logist11}, an electrical circuit \cite{Teleke10}, the power split in a hybrid electric vehicle \cite{Sciarretta04}, or a biological process \cite{Sharp21}.

A well-known and general way to solve optimal control problems is by using Bellman's dynamic programming \cite{Bellman:1957} which in the continuous-time case reduces to solving the value function $V(x,t)$ from the Hamilton--Jacobi--Bellman (HJB) equation (see \eqref{eq:HJB_equation} below) \cite{Lewis_book12,Kirk:2004}. The value function then completely determines the optimal feedback law (see \eqref{eq:optimal_control} below). However, except for the closed form solution to the linear quadratic tracker (LQT) case (see Sec.~\ref{sub:LQT_control_problem}), solving the HJB equation is usually computationally heavy \cite{Tahirovic22}, and in most cases, inherently an offline procedure. Computational methods for solving the HJB equation include, for example, finite-difference methods \cite{Crandall84}, reachable set methods \cite{Bokanowski13}, neural networks \cite{Abu-Khalaf05,Kim20}, and the successive approximation approach \cite{Beard98,Maruta20}. 

\begin{figure}[tb]
\centering
\includegraphics[width=\columnwidth]{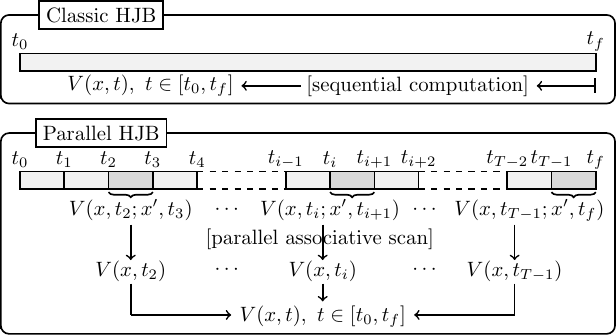}
\caption{Illustration of the basic idea. Instead of obtaining the value function from the HJB equation on a (sequential) backward pass, we first solve the conditional value functions on sub-intervals and then combine them in parallel using a parallel associative scan.}
\label{fig:concept}
\end{figure}

In this paper, the aim is to develop parallel-in-time methods for solving the HJB equation which can be used to speedup the computations in highly parallel architectures such as multi-core central processing units (CPUs) or graphics processing units (GPUs). The basic idea of the present approach is illustrated in Fig.~\ref{fig:concept}. We first split the time-interval $t \in [t_0,t_f]$ into $T$ sub-intervals $t_0 < t_1 < \cdots < t_T = t_f$. We then solve the conditional value functions $V(x,t_i;x',t_{i+1})$ at each sub-interval, which can be done fully in parallel in $O(1)$ span time complexity. These conditional value functions then obey an associative combination rule which enables the parallel computation of the (non-conditional) value functions at any given time on the interval $[t_0,t_f]$. Due to the associative property, this step can be implemented by using parallel scan algorithms \cite{Ladner:1980,Blelloch:1989,Blelloch:1990} in $O(\log T)$ span time complexity. The main methodological challenges in the present approach are related to representation and computation of the conditional value functions along with the implementation of the combination rule, which are thus the foci of this paper. 

Our approach for parallelisation uses a similar framework of parallel associative scans as was used for discrete-time dynamic programming based optimal control in \cite{Sarkka:2022}, which in turn is related to the parallel state-estimation methodology in \cite{Sarkka:2021, Hassan:2021, yaghoobi2021parallel, Corenflos:2022}. However, as opposed to these previous works, in this paper, we consider continuous-time problems. 

For time-parallelisation of dynamic programming (i.e., HJB) for continuous-time problems, an alternative approach would be to use parallel-in-time differential equation solving methods \cite{Gander:2015} such as parareal \cite{Lions:2001,gander2007analysis} for HJB or LQT equations. In particular, parareal \cite{Lions:2001,gander2007analysis} is a state-of-the-art method which is based on splitting the interval into sub-intervals and then iterating coarse and dense solvers such that the result converges to the dense solver result. Although the dense solution of parareal is fully parallelisable, the coarse solver is not and hence the span time complexity of parareal is still $O(T)$. Many HJB solving methods such as finite-differences also allow for spatial-direction parallelisation of individual time steps in the algorithm, but this still keeps the time complexity in $O(T)$. This kind of parallelisation could also used within the proposed methodology.

Parareal has also been previously used to parallelise the solutions of forward and backward differential equations arising from classical Pontryagin/Euler--Lagrange type of solutions to continuous-time optimal control problems. In particular, \cite{Maday+Turinici:2002,Maday:2007,Riahi:2016} apply this kind of method to quantum optimal control problems and \cite{maday2002parareal,Gander:2020} to optimal control of ordinary and partial differential equations. However, as opposed to these references, the focus of the present paper is in parallel HJB-based methods for optimal control problems.

The main contributions of this paper are:
\begin{itemize}
\item Framework for span $O(\log T)$ time-parallelisation of the HJB equation solution via an associative combination of conditional value functions.
\item Parallel algorithm to solve LQT problems in span $O(\log T)$ time by specialising the aforementioned approach to LQT problems.
\item Analytical and numerical solution methods, and experimental evaluation of the performance of the methods on a multi-core CPU and a GPU.
\end{itemize}

The structure of the paper is the following. In Section~\ref{sec:background} we briefly review the HJB and LQT solutions to optimal control problems, and parallel all-prefix-sums computing methods using associative scans. Section~\ref{sec:parallel_formulation} introduces the associative formulation of general optimal control problems and presents the backward and forward equations for the conditional value functions. In Section~\ref{sec:lqt} we specialise the general associative formulation to LQT problems and derive the backward and forward differential equations for the value function of LQT. Section~\ref{sec:direct_opt} discusses the computation of conditional value functions via direct optimisation. In Section~\ref{sec:results} we experimentally evaluate the speed benefit of the parallelisation on multi-core CPU and GPU. Section~\ref{sec:conclusion} finally concludes the paper.

\section{Background} \label{sec:background}

This section reviews the required background in optimal control and parallel computing. The continuous-time control problem and the LQT case are reviewed in Sections \ref{sub:control_problem} and \ref{sub:LQT_control_problem}, respectively. The required background on parallel computing is provided in Section \ref{sub:assoc_operators}.

\subsection{Continuous-time optimal control problem}\label{sub:control_problem}

We consider continuous-time optimal control problems on a fixed time interval $t \in [t_0,t_f]$. We denote the system state and control at time $t$ as $x(t)\in \mathbb{R}^{n_x}$ and $u(t)\in \mathbb{R}^{n_u}$, respectively. The ordinary differential equation (ODE) characterising the dynamics, and the cost functional to be minimised are respectively given as (cf.\ \cite{Lewis_book12})
\begin{equation}
\label{eq:Control_problem}
\begin{split}
    \frac{dx(t)}{dt} &= f(x(t), u(t), t), \\
  C[u] &= \phi(x(t_f)) + \int_{t_0}^{t_f} \ell(x(t),u(t),t) \, dt,
\end{split}
\end{equation}
where $f(\cdot, \cdot, t)$ is the function that models the system dynamics at time $t$, $\ell(\cdot,\cdot,t)$ is the weighting function at time $t$, and $\phi(x(t_f))$ is the final-time weighting function. %

The optimal control problem consists of finding a control $u^*(t)$ for $t \in [t_0,t_f]$ to drive the system, with a known initial state $x(t_0)$, along a trajectory $x^*(t)$ such that the cost $C[u]$ is minimised. We assume that the optimal control problem \eqref{eq:Control_problem} is well-posed in the sense that there exists a time-measurable control function $u^*(t)$ which minimises the cost-functional. However, this solution does not need to be unique. There are various assumptions that can be used to ensure that an optimal control exists \cite{Fleming+Rishel:1975,Bardi_book97,Fleming_book06}. Typically, we assume that $f$ is locally Lipschitz in $x$-variable, uniformly in $u$, and continuous with respect to $t$. Additionally, $\ell$ and $\phi$ need to be bounded from below, continuously differentiable in $x$ and $u$, and continuous in time. 

In this paper, we are interested in finding an optimal feedback control $u^*(t) = u^*(x^*(t),t)$ in which case it is natural to formulate the problem in terms of Bellman's dynamic programming \cite{Bellman:1957,Lewis_book12,Kirk:2004}. In this approach, the optimal cost $V(x,t)$ of being at state $x$ at time $t$, also called the value function, is obtained by solving the HJB equation \cite{Lewis_book12}
\begin{equation}
\label{eq:HJB_equation}
\begin{split}
-\frac{\partial V(x,t)}{\partial t}
&= \min_{u} \left\{
\ell(x,u,t)
+ \left( \frac{\partial V(x,t)}{\partial x} \right)^\top \, f(x,u,t)
\right\}
\end{split}
\end{equation}
backwards in time from the boundary condition $V(x,t_f) =  \phi(x)$. Then, the optimal feedback control is \cite{Lewis_book12}
\begin{equation}
\label{eq:optimal_control}
\begin{split}
  u^*(x,t) = \arg \min_{u} \left\{
    \ell(x,u,t)
    + \left( \frac{\partial V(x,t)}{\partial x} \right)^\top \, f(x,u,t)
    \right\}.
\end{split}
\end{equation}

The HJB equation \eqref{eq:HJB_equation} formally requires that the value function $V(x,t)$ is continuously differentiable in $x$. In many cases it is not, and therefore one needs to allow for so-called viscosity solutions \cite{Bardi_book97,Fleming_book06} to the equation. Viscosity solutions are allowed to have discontinuous derivatives in some parts of the state space. The existence conditions and properties of viscosity solutions can be found in \cite{Bardi_book97,Fleming_book06} and references therein, and they are more restrictive than the mere existence of an optimal control solution. In this paper, the solution to the HJB equation is also understood in viscosity sense although in some occasions of the paper, which we try to point out explicitly at the text, the differentiability is assumed.

Given an optimal control, the optimal trajectory is obtained by solving the ODE
\begin{equation}
\label{eq:optimal_trajectory}
\begin{split}
\frac{dx^*(t)}{dt} &= f^*(x^*(t), t),\\
&= f(x^*(t), u^*(x^*(t),t), t)
\end{split}
\end{equation}
starting at $x^*(t_0)=x(t_0)$. We also assume that this solution exists, which can be ensured by, for example, requiring that function $x \mapsto f(x, u^*(x,t), t)$ is locally Lipschitz and $t \mapsto f(x, u^*(x,t), t)$ is continuous.

In general, the HJB equation \eqref{eq:HJB_equation} and the optimal control \eqref{eq:optimal_control} do not have closed-form expressions and they need to be approximated. One type of solution is to discretise the time interval and state space, and apply a finite difference approximation \cite{Crandall84} or related methods \cite{Osher91,Zhang03,Tsitsiklis95}. Then, for the optimal controls at the discretised time steps, one can use a numerical ODE solver \cite{Hairer_book08} to compute the optimal trajectory \eqref{eq:optimal_trajectory}. In the next section, we review the LQT problem, for which we can obtain a solution in terms of ODEs, without discretising the state space.

\subsection{Continuous-time LQT problem}\label{sub:LQT_control_problem}
In the LQT problem, the system \eqref{eq:Control_problem} has the form \cite{Anderson:1989,Lewis_book12}
\begin{equation}
\label{eq:LQT_Control_problem}
\begin{split}
f(x,u,t) &= F(t) \, x + L(t) \, u + c(t), \\
\ell(x,u,t) &= \frac{1}{2} (r(t) - H(t) \, x)^\top \, X(t) \, (r(t) - H(t) \, x)\\
&+ \frac{1}{2} u^\top U(t) \, u, \\
\phi(x) &= \frac{1}{2} (H_f \, x - r_f)^\top \, X_f \, (H_f \, x - r_f), \\
\end{split}
\end{equation}
where for $t \in [0,t_f]$, $r(t) \in \mathbb{R}^{n_r}$ is the reference trajectory, $r_f \in \mathbb{R}^{n_r}$ the target point, $H(t)$ is an $n_r \times n_x$ output matrix, $H_f$ is an $n_r \times n_x$ final output matrix, $F(t)$ is an $n_x \times n_x$ matrix, $L(t)$ an $n_x \times n_x$ matrix, and $c(t)$ is a $n_x$-dimensional vector. In addition, $X(t)$ and $X_f$ are symmetric positive semi-definite $n_r \times n_r$ matrices and $U(t)$ is an $n_u \times n_u$ positive definite symmetric matrix. We also assume that $r(t)$, $H(t)$, $F(t)$, $L(t)$, $c(t)$, $X(t)$, and $U(t)$ are continuous functions of time. With these conditions, the unique solution to the optimal control problem is ensured and the value function is a continuously differentiable non-negative function (see, e.g., \cite{Anderson:1989,Lewis_book12}). Please note that the controllability of the pair $(F,L)$ is not required, because the cost will be finite provided that the control is bounded and the (strict) positive definiteness of $U(t)$ ensures that a unique minimum always exists.

The objective in \eqref{eq:LQT_Control_problem} can be interpreted such that a linear combination $H(t) \, x(t)$ of the trajectory $x(t)$ follows the reference trajectory $r(t)$. The linear quadratic regulator (LQR) problem \cite{Lewis_book12} can be recovered by setting $c(t)=0$, $r(t) = r_f = 0$, and $H(t) = H_f = I$, where $I$ is an identity matrix of size $n_x \times n_x$.  

The HJB equation \eqref{eq:HJB_equation} for the LQT problem \eqref{eq:LQT_Control_problem} is solved by an optimal cost function $V(x,t)$ of the form \cite{Lewis_book12}
\begin{equation}\label{eq:V_quadratic_LQT}
\begin{split}
V(x,t) = \zeta(t) + \frac{1}{2} x^\top \, S(t) \, x - v^\top(t) \, x,
\end{split}
\end{equation}
where $v(t)$ is an $n_x$-dimensional vector, $S(t)$ is an $n_x \times n_x$ symmetric matrix, and $\zeta(t)$ is a function that does not depend on the state and has no effect on the optimal control, so it does not need to be computed. 

The parameters $v(t)$ and $S(t)$ in \eqref{eq:V_quadratic_LQT} can be solved backwards, starting at $t_f$, by first setting $v(t_f) = H_f^\top \, X_f \, r_f$, $S(t_f) = H_f^\top \, X_f \, H_f$ and then solving following ODEs backwards until time $t = t_0$:
\begin{equation}
\label{eq:LQT_differential_S_v}
\begin{split}
\frac{dS(t)}{dt} &=
- F^\top(t) \, S(t)
- S(t) \, F(t) - H^\top(t) \, X(t) \, H(t)
 \\
& \quad  + S(t) \, L(t) \, U^{-1}(t) \, L^\top(t) \, S(t),  \\
\frac{dv(t)}{dt}  &=
- H^\top(t) \, X(t) \, r(t)
+ S(t) \, c(t) 
- F^\top(t) \, v(t) \\
& \quad + S(t) \, L(t) \, U^{-1}(t) \, L^\top(t) \, v(t).
\end{split}
\end{equation}
Given $v(t)$ and $S(t)$, the optimal control law \eqref{eq:optimal_control} becomes
\begin{equation}
\label{eq:LQT_optimal_control}
\begin{split}
	u^*(x,t) = - U^{-1}(t) \, L^\top(t) \, S(t) \, x + U^{-1}(t) \, L^\top(t) \, v(t).
\end{split}
\end{equation}
Substituting \eqref{eq:LQT_optimal_control} and $f(x,u,t)$ in \eqref{eq:LQT_Control_problem} into \eqref{eq:optimal_trajectory}, the ODE for the optimal trajectory becomes
\begin{equation}
\label{eq:optimal_trajectory_ODE_LQT}
\begin{split}
\frac{dx^*(t)}{dt}&= \widetilde{F}(t) \, x^*(t) + \widetilde{c}(t), \\
\widetilde{F}(t)&= F(t) - L(t) \, U^{-1}(t) \, L^\top(t) \, S(t),   \\  
\widetilde{c}(t)&= L(t) \, U^{-1}(t) \, L^\top(t) \, v(t) + c(t),
\end{split}
\end{equation}
with the initial condition $x^*(t_0)=x(t_0)$.

The solution to \eqref{eq:optimal_trajectory_ODE_LQT} can be written in terms of an $n_x \times n_x$ transition matrix $\Psi(t,t_0)$ and $n_x$-dimensional vector $\alpha\left(t,t_{0}\right)$ such that (see, e.g., Sec.\ 2.3 of \cite{Sarkka_book19})
\begin{equation}\label{eq:optimal_trajectory_LQT}
x^{*}(t)=\Psi(t,t_{0}) \, x(t_{0})+ \alpha\left(t,t_{0}\right),
\end{equation}
where
\begin{align}
\frac{\partial\Psi(t,t_{0})}{\partial t} &= \widetilde{F}(t) \, \Psi(t,t_{0}), \label{eq:ODE_transition_matrix} \\ 
\alpha\left(t,t_{0}\right) &= \int_{t_{0}}^{t} \Psi(t,\tau) \, \widetilde{c}(\tau) \, d\tau. \label{eq:ODE_transition_matrix2}
\end{align}

Applying the Leibniz integral rule, we can write \eqref{eq:ODE_transition_matrix2} as the ODE
\begin{equation}\label{eq:ODE_phi}
\frac{\partial\alpha\left(t,t_{0}\right)}{\partial t}=\widetilde{F}(t) \, \alpha\left(t,t_{0}\right) + \widetilde{c}(t).
\end{equation}
We can then approximate $\Psi(t,t_{0})$ and $\alpha\left(t,t_{0}\right)$ by using an ODE solver for \eqref{eq:ODE_transition_matrix} and \eqref{eq:ODE_phi} with initial conditions $\Psi(t_{0},t_{0})=I$ and $\alpha\left(t_{0},t_{0}\right)=0$. Note that we usually only know $S(t)$ and $v(t)$ at some discretised time steps, as they are solutions to ODEs that we solve numerically. Therefore, the ODEs \eqref{eq:ODE_transition_matrix} and \eqref{eq:ODE_phi} are typically solved by discretising at the same time steps.

\subsection{Associative operators and parallel computing}\label{sub:assoc_operators}
In this paper we use parallel scans, also called all-prefix-sums or scan algorithms \cite{Blelloch:1989,Blelloch:1990} to parallelise the continuous-time control problem. We proceed to review the main ideas behind parallel scans. Given a sequence of elements $a_1,\ldots,a_T$ and an associative operator $\otimes$ on them, the all-prefix-sums operation or scan computes $s_1,\ldots, s_T$ such that
\begin{equation}
\begin{split}´
s_1 &= a_1, \\
s_2 &= a_1 \otimes a_2, \\
&~ \vdots \\
s_T &= a_1 \otimes a_2 \otimes \cdots \otimes a_T.
\end{split}
\end{equation}
The associative operator can be a sum, a multiplication, a maximisation, or any other associative operator. Being associative means that the operator satisfies $(a \otimes b) \otimes c = a \otimes (b \otimes c)$ for any three elements $a,b,c$ for which the operator is defined. We can also similarly compute reverse all-prefix-sums $\bar{s}_k = a_k \otimes a_{k+1} \otimes \ldots \otimes a_T$. For example, when the sequence of elements is
\begin{equation}
\begin{split}
   a_{1:T} = [ 1, 2, 3, 4, 5, 6, 7, 8 ],
\end{split}
\label{eq:a_1T}
\end{equation}
where thus $T = 8$, and the associative operator is summation $\otimes = +$, then the all-prefix-sums is the cumulative sum
\begin{equation}
\begin{split}
   s_{1:T} = [ 1,     3,     6,    10,    15,    21,    28,    36 ],
\end{split}
\label{eq:s_1T}
\end{equation}
and the reverse all-prefix-sums is
\begin{equation}
\begin{split}
   \bar{s}_{1:T} = [ 36,    35,    33,    30,    26,    21,    15,     8 ].
\end{split}
\end{equation}
A direct sequential way to compute the all-prefix-sums in \eqref{eq:s_1T} is by the loop:
\begin{itemize}
\item $s_0 = 0$ (or the neural element of the general operator $\otimes$).
\item For $i = 1,\ldots,T$ do
\begin{itemize}
\item[] $s_{i} = s_{i-1} \otimes a_i$.
\end{itemize}
\end{itemize}
However, the direct sequential computation of the all-prefix-sums operation inherently takes $O(T)$ time, even when run on a parallel computer. The parallel scan algorithm \cite{Blelloch:1989,Blelloch:1990} is an algorithm that uses parallelisation to compute the all-prefix-sums operation in $O(\log T)$ span time, which corresponds to $O(\log T)$ parallel operation steps. To do so, it performs up-sweep and a down-sweep computations on a binary tree that generates independent sub-problems that can be computed in parallel. The specific algorithm is, for example, given in \cite[Alg. 1]{Sarkka:2022}. We proceed to illustrate how parallel scans work in the following example.

\begin{figure}[tb]
	\centering
	\includegraphics[width=\columnwidth]{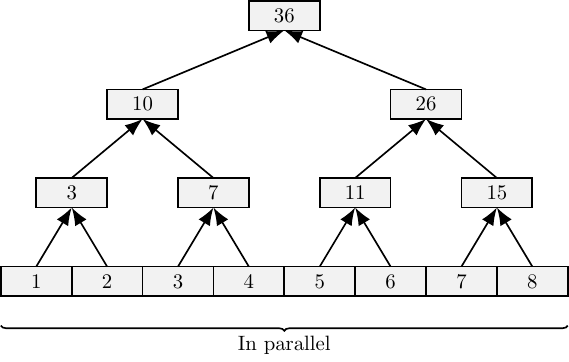}
	\caption{Illustration of the up-sweep of the parallel scan algorithm to compute the all-prefix-sums of $a_{1:T} = [ 1, 2, 3, 4, 5, 6, 7, 8 ]$.}
	\label{fig:parallel_up}
\end{figure}

\begin{figure}[tb]
	\centering
	\includegraphics[width=\columnwidth]{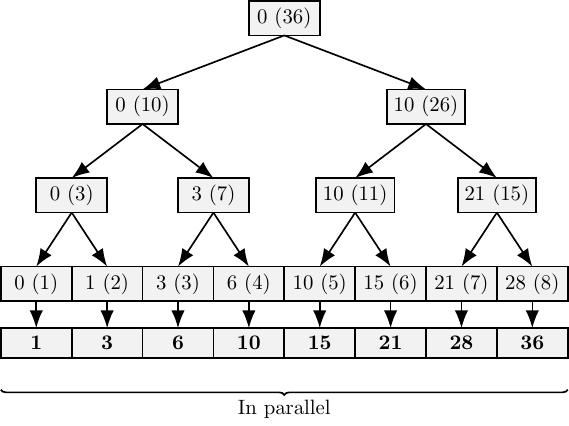}
	\caption{Illustration of the down-sweep of the parallel scan algorithm to compute the all-prefix-sums of $a_{1:T} = [ 1, 2, 3, 4, 5, 6, 7, 8 ]$. The results of the up-sweep algorithm are shown in parentheses, see Figure \ref{fig:parallel_up}. The bottom row shows the computation of the all-prefix-sums.}
	\label{fig:parallel_down}
\end{figure}

\begin{exa}
	Consider the sequence of elements given in \eqref{eq:s_1T} with the operator $+$. The direct sequential computation now takes $T = 8$ steps where we apply operator $\otimes$ once at each step. The operation of the up-sweep of the parallel scan algorithm is illustrated in Fig.~\ref{fig:parallel_up}. We start from the bottom of the tree by computing the sums of pairs, then the sums of pairs of the results, and finally the total sum to the root element. Each of the levels of tree can be processed in parallel which requires total of $\log_2 8 = 3$ steps. The down-sweep is illustrated in Fig.~\ref{fig:parallel_down}. In that algorithm, we assign a value $L$ to each of the nodes as follows:
\begin{itemize}
\item The root has $L = 0$ (or the neural element of $\otimes$).
\item Every left child inherits the present value $L$.
\item Every right child gets the left child's up-pass sum plus the current $L$.
\item At each leaf, we output $s_i = L + a_i$.
\end{itemize}
This down-sweep again takes $\log_2 8 = 3$ steps plus the leaf operations which sums to $7$ steps. Although in this example, the difference between the number of steps in the sequential and tree versions is not significant (or is even none if we start the loop from $s_1 = a_1$), in general, the up and down sweeps take total of $2\log T + 1$ steps which is significantly less than $T$ when $T$ is large.
\end{exa}

\section{Parallel formulation of continuous-time control problems} \label{sec:parallel_formulation}

This section presents a parallel formulation of the continuous-time control problem \eqref{eq:Control_problem} based on associative operators. Section \ref{subsec:Conditional_value_function} defines conditional value functions, and Section \ref{subsec:Cond_init} their initialization. Section \ref{subsec:Cond_back_forw_HJB} presents conditional backward and forward HJB equations for the conditional value functions. Section \ref{subsec:Parallel_opt_control} explains how to solve optimal control problems using parallel scans. Optimal trajectory recovery is addressed in Section \ref{subsec:Optimal_trajectory}.

\subsection{Conditional value functions and combination rule}\label{subsec:Conditional_value_function}
We introduce the conditional value function via the following definition. With a slight abuse of notation, both value functions and conditional value functions are denoted by letter $V$, though value functions have two arguments and conditional value functions have four arguments. Let $\mathcal{R}(x,s;\tau)$ denote the set of states that can be reached at time $\tau$ from the state $x(s)$ using an admissible time-measurable control function $u(t)$.

\begin{defn}[Conditional value function]\label{Defn:conditional_value_function}
	For $x\in\mathbb{R}^{n_{x}}$ and $z\in\mathbb{R}^{n_{x}}$, the conditional value function $V(x,s;z,\tau)$ from time $s \in [t_0,t_f]$ to time $\tau \in [t_0,t_f]$ is the cost of the optimal trajectory in time interval $[s,\tau]$ starting at state $x(s)=x$ and ending at state $x(\tau)=z$. When $z \in \mathcal{R}(x,s;\tau)$, then we have 
	\begin{equation}
	\label{eq:conditional_value_opt}
	\begin{split}
	V(x,s;z,\tau) &= \min_{u} \left\{
	\int_s^\tau \ell(x(t),u(t),t) \, dt 
	\mid x(\tau) = z \right\} \\
	\end{split}
	\end{equation}
subject to %
\begin{equation}
\label{eq:conditional_value_con}
\begin{split}
    \frac{dx(t)}{dt} &= f(x(t), u(t), t), 
\end{split}
\end{equation}
on interval $[s,\tau]$ with given $x(s) = x$, where the notation $|$ in the minimum indicates a constraint $x(\tau) = z$. When $z \notin \mathcal{R}(x,s;\tau)$, then we set $V(x,s;z,\tau) = \infty$.
\end{defn}

Given two conditional value functions, one from time $s$ to time $\tau$, and one from time $\tau$ to time $t$, we can combine them to get a conditional value function from time $s$ to $t$. To achieve this, we first define the combination operation as follows. 
\begin{defn}\label{Def:Assoc_operator_HJB}
Given functions $W_{1}, W_{2}: \mathbb{R}^{n_{x}} \times \mathbb{R}^{n_{x}} \rightarrow \mathbb{R}$, we define
\begin{equation}
W_{1}(x, \cdot) \otimes W_{2}(\cdot, y):=\min _{z}\left(W_{1}(x, z)+W_{2}(z, y)\right).
\end{equation}
In addition, given a function $W: \mathbb{R}^{n_{x}} \times \mathbb{R}^{n_{x}} \rightarrow \mathbb{R}$ and a function $\varphi: \mathbb{R}^{n_{x}} \rightarrow \mathbb{R}$, we define
\begin{equation}
W(x, \cdot) \otimes \varphi:=\min _{z}(W(x, z)+\varphi(z)).
\end{equation}
\end{defn}

This combination operation is associative \cite{Sarkka:2022}, as required in the parallel scan algorithm, see Section \ref{sub:assoc_operators}. The calculation of the conditional value functions based on this operator is given in the following theorem.

\begin{thm} \label{the:HJB-comb}
	Given $s<\tau<t$, and for any $x, y \in \mathbb{R}^{n_{x}}$, the conditional value function from time $s$ to time $t$ satisfies
	\begin{equation}\label{eq:conditional_V_operator}
	\begin{split}
	V(x, s ; y, t) &=V(x, s ; \cdot, \tau) \otimes V(\cdot, \tau ; y, t). \\
	\end{split}
	\end{equation}
	Furthermore we have
	\begin{equation}\label{eq:V_operator}
	 V(x, s)=V(x, s ; \cdot, \tau) \otimes V(\cdot, \tau) .
	\end{equation}
\end{thm}

The combination rules given in Theorem~\ref{the:HJB-comb} are discrete-time in the sense that we need to know the conditional value functions, such as $V(x,s;z,\tau)$, already for some finite, non-vanishing time-interval $\tau - s > 0$. The computation of these quantities is considered in Sections \ref{subsec:Cond_init} and \ref{subsec:Cond_back_forw_HJB}.

\begin{rem}
	It should be noted that Definition \ref{Defn:conditional_value_function} only considers the case of a fixed end-point and no final weighting function. However, it is possible to include the final weighting function to the current methodology by introducing a time $t_f^+$ which is ``infinitesimally after'' $t_f$ and for which the conditional value function satisfies
	\begin{equation}
	V(x,t_f;z,t_f^+)=V(x,t_f)=\phi(x).
	\end{equation}	
\end{rem}

\subsection{Initialization of conditional value functions} \label{subsec:Cond_init}

To implement the combination rule in the previous section, we need to compute the conditional value function for a finite time interval $\tau - s > 0$. We can now notice that the problem defining the conditional value function in \eqref{eq:conditional_value_opt} and \eqref{eq:conditional_value_con} is a problem of the form \eqref{eq:Control_problem} on an interval $[s,\tau]$ with the initial value $x$ at time $s$, the terminal value $z$ at time $\tau$, and with a final-time weighting function $\phi(x) = 0$. The conditional value function is simply the optimal cost for this problem.

Thus for each pair $(x,z)$, the optimal trajectory $x^*(t)$ and the optimal control $u^*(t)$ for the aforementioned problem can computed via classical means by augmenting the dynamic model and the terminal conditions to the cost function via Lagrange multipliers and solving the minimisation problem via calculus of variations (see, e.g., \cite{Lewis_book12,Bryson+Ho:1975}). The conditional value function is given by the cost evaluated at the optimal trajectory and control:
\begin{equation}
\begin{split}
  V(x,s;z,\tau) = \int_s^\tau \ell(x^*(t),u^*(t),t) \, dt,
\end{split}
\end{equation}
where the dependence on $x$ and $z$ comes through the initial and terminal conditions on $x^*(t)$. 

The solution to this problem can also be obtained numerically, for example, by using direct or indirect shooting methods \cite{Bock84,Betts98}. We discuss such numerical approaches in Section~\ref{sec:direct_opt}.

\subsection{Backward and forward conditional HJB equations} \label{subsec:Cond_back_forw_HJB}
As discussed in the previous section, for a given $s$ and $\tau$ we can, in principle, solve the conditional value function in Definition \ref{Defn:conditional_value_function}, for example, by using calculus of variations. However, we can also form backward and forward HJBs for the conditional value functions as will be seen in this section.

Let us now assume that we have solved the conditional value function for a small time interval $[s',\tau]$, for example, using the variational calculus. Provided that the conditional value function $V(x,s';z,\tau)$ is differentiable in $x$ and $s'$, then by taking $x \to V(x,s';z,\tau)$ as the final weighting function at $s'$, we get that for $s \le s'$ the conditional value function $V(x,s;z,\tau)$ solves the backward conditional HJB equation
\begin{equation}\label{eq:backward_conditional_HJB}
\begin{split}
&-\frac{\partial V(x,s;z,\tau)}{\partial s}\\
&= \min_{u} \left\{
\ell(x,u,s)
+ \left( \frac{\partial V(x,s;z,\tau)}{\partial x} \right)^\top f(x,u,s)
\right\}.
\end{split}
\end{equation}
This equation goes backwards in time reducing the lower time limit $s$ in $V(x,s;z,\tau)$ and it is defined only for arguments $x,s;z,\tau$ such that $z \in \mathcal{R}(x,s;\tau)$, that is, $V(x,s;z,\tau) < \infty$. When $V(x,s;z,\tau)$ is not differentiable in $x$ everywhere, we can still interpret this solution in a viscosity sense \cite{Bardi_book97,Fleming_book06}.

When $V(x,s;z,\tau)$ is differentiable in $z$ and $\tau$, as shown in Appendix \ref{sec:Proof_forward_HJB_append}, at least in a viscosity sense, we can also form the corresponding forward conditional HJB equation to increase the upper time limit $\tau$:
\begin{equation}\label{eq:forward_conditional_HJB}
\begin{split}
&\frac{\partial V(x,s;z,\tau)}{\partial \tau}\\
&= \min_{u} \left\{
\ell(z,u,\tau)
- \left( \frac{\partial V(x,s;z,\tau)}{\partial z} \right)^\top f(z,u,\tau)
\right\}.
\end{split}
\end{equation}

These results suggest a computational procedure where we first solve the conditional value function for a small time interval $[s',\tau]$ using the calculus of variations (or its numerical approximation), then proceed with one of the HJBs above with the small time interval solution as the boundary condition, and finally use \eqref{eq:V_operator} to incorporate the boundary condition at the end time. However, we can also consider what happens to the conditional value function when $s' \to \tau$. In this limit, the cost of the path becomes zero, but at the same time the only end point $z$ that we can reach from an initial point $x$ is the one with $z = x$. Thus, at least formally, we have
\begin{equation}
\label{eq:V_zero_length}
V(x,\tau;z,\tau) =\begin{cases}
0 & x=z\\
\infty & x\neq z
\end{cases},
\end{equation}
where the $\infty$ results from our convention of assigning infinite value to the state pairs that are not reachable from each other. Thus, at least in principle, we can directly use this as the boundary condition for the HJBs above without a need to compute an initial solution for a small time interval. However, handling the infinities can be numerically challenging and the solution to the HJB might not exist in any meaningful sense.

Another way to interpret the above boundary condition is to see it as the penalty (or barrier) function to implement the terminal constraint, which corresponds to
\begin{equation*}
\begin{split}
&V(x,s;z,\tau) \\
&= \min_{u} \Bigg\{ \frac{1}{2\epsilon} (x(\tau) - z)^\top (x(\tau) - z) %
+ 
\int_s^\tau \ell(x(t),u(t),t) \, dt
\Bigg\},
\end{split}
\end{equation*}
with $\epsilon \to 0$. Because in the limit, all the trajectories with $x(\tau) \ne z$ will have an infinite cost, the optimal cost (which is going to be finite) has to be associated with a trajectory which has $x(\tau) = z$.

It is easy to check (cf. \cite{Dorato+Abdallah+Cerone:1995}) that in the LQT case, the above indeed gives the terminal-constrained solution. Therefore, we also use this formulation in the parallel LQT presented in Section~\ref{sec:lqt}. However, in general, the penalty function formulation might also cause the HJB solution not to exist.

\subsection{Parallel solution of optimal control}\label{subsec:Parallel_opt_control}

To solve the optimal control in parallel, we first split the time interval $[t_0,t_f]$ into $T$ sub-intervals $[t_0,t_1],\ldots,[t_{T-1},t_T]$, where $t_T=t_f$. We then define an associative element $a_j$ in each interval and recover the value functions via the following theorem.

\begin{thm} \label{the:parallel_HJB}
If we let the associative elements $a_j$, for $j\in \{0,\ldots,T \}$, be the conditional and final value functions as follows:
\begin{align}
a_{j}&=V(\cdot,t_{j};\cdot,t_{j+1}), \quad j\in\{0,\ldots,T-1\}, \label{eq:a_j}\\
a_{T}&=V(\cdot,t_{T}),
\end{align}
then
\begin{align}
a_0 \otimes a_1 \otimes \cdots \otimes a_k &= V(\cdot,t_0;\cdot,t_{k+1}), \label{eq:parallel_cond_B_HJB} \\ 
a_k \otimes a_{k+1} \otimes \cdots \otimes a_T &= V(\cdot,t_{k}), \label{eq:parallel_V_HJB}
\end{align}
where the associative operator $\otimes$ is defined in Theorem~\ref{the:HJB-comb}.
\end{thm}

This theorem is a direct consequence of applying Theorem~\ref{the:HJB-comb} recursively. The element initiation can be done in parallel for each $j$. That is, element $a_T$ is given by $V(\cdot,t_{T})=\phi(\cdot)$, and the elements $a_{j}$ for $j<T$ can be obtained by using the classical calculus of variations solution or by solving either the backward or forward HJB equations in the interval $[t_j,t_{j+1}]$. Then, based on \eqref{eq:parallel_cond_B_HJB}, \eqref{eq:parallel_V_HJB} and the associative operator in Definition \ref{Def:Assoc_operator_HJB}, we can use parallel associative scans to obtain $V(\cdot,t_0;\cdot,t_k)$ and $V(\cdot,t_{k})$. Finally, we can use the forward or backward HJB equations to find the value function at any given time within each sub-interval, which can be done in parallel.

After computing  $V(\cdot,t)$ for all $t \in  [t_0,t_f]$, we can compute the optimal control at each time using \eqref{eq:optimal_control}. This can be done in parallel for all times.

\subsection{Optimal trajectory recovery}\label{subsec:Optimal_trajectory}

In the previous section we computed the value functions and the optimal control law in parallel using parallel associative scans. We can now proceed to compute the optimal trajectory $x^*(t)$ which is the solution to ODE \eqref{eq:optimal_trajectory}. As in discrete time control problems \cite{Sarkka:2022}, there are two methods for the optimal trajectory recovery.

\subsubsection{Method 1}
In parallel, for each sub-interval $j$, we form $\psi^*(x(t_{j-1}),t)$, which is the solution of the ODE \eqref{eq:optimal_trajectory} for $t \in [t_{j-1},t_j]$ with an initial state $x(t_{j-1})$. For notational simplicity, we denote by $\psi^*_j(x(t_{j-1})) = \psi^*(x(t_{j-1}),t_j)$ the resulting function that provides the optimal state at time $t_{j}$ given the state at time $t_{j-1}$.

The optimal trajectory $x^*(t_j)$ at times $t_j$ for $j=0,\ldots,T$ can be calculated by function composition
\begin{equation}
x^*(t_j) = (\psi^*_j \circ \psi^*_{j-1} \circ \ldots \circ \psi^*_1)(x(t_0)).
\end{equation}
Because function composition is associative, calculating $x^*(t_j)$ for $j=0,\ldots,T$ can be done using parallel scans \cite[Sec. III.C]{Sarkka:2022}. The  optimal trajectory for a time $t \in [t_{j-1},t_j]$ can be recovered using $\psi^*(x^*(t_{j-1}),t)$, which can be done in parallel in each sub-interval.

\subsubsection{Method 2}
From the definition of value function and conditional value function, we directly obtain that the value of the optimal trajectory at time $t$ is
\begin{equation}\label{eq:Opt_trajectory_method2}
x^{*}(t)=\arg\min_{x}\{V(x(t_{0}),t_{0},x,t)+V(x,t)\}.
\end{equation}
This equation directly provides the means to compute $x^{*}(t)$ in parallel for all $t$.

\section{Parallelisation of continuous-time LQT} \label{sec:lqt}

This section explains the parallelisation of the continuous-time LQT problem, by making use of the general results in the previous section. The temporal parallelisation of the LQT is an important result as it can be addressed solving ODEs in combination with parallel scans. In the LQT problem, the conditional value functions and combination rule are the same as in the previous section, though they admit specific parameterisations. These are provided in Section \ref{sub:conditional_value_LQT}. The backward and forward conditional HJB equations are explained in Sections \ref{sub:conditional_backward_LQT} and \ref{sub:conditional_forward_LQT}, respectively. The parallel solution of continuous-time LQT is addressed in Section \ref{sub:parallel_LQT}. Optimal trajectory recovery is provided Section \ref{sub:optimal_trajectory_LQT}.

\subsection{Conditional value functions and combination rule}\label{sub:conditional_value_LQT}

For the LQT problem \eqref{eq:LQT_Control_problem}, the solution to the HJB equation is of the form \eqref{eq:V_quadratic_LQT}. It will be proved in Sections \ref{sub:conditional_backward_LQT} and \ref{sub:conditional_forward_LQT} that the conditional value function $V(x,s;z,\tau)$ for the LQT problem is the solution to the optimisation problem
\begin{equation}\label{eq:LQT_V_form1}
V(x,s;z,\tau)=\max_{\lambda}g\left(x,s;z,\tau,\lambda\right),
\end{equation}
where $\lambda$ is an $n_x$-dimensional vector with the Lagrange multipliers and $g\left(x,s;z,\tau,\lambda\right)$ is the dual function, which has the form 
\begin{align}
&g\left(x,s;z,\tau,\lambda\right) =\zeta(s,\tau)+\frac{1}{2}x^{\top}J\left(s,\tau\right)x-x^{\top}\eta\left(s,\tau\right) \nonumber \\
&-\frac{1}{2}\lambda{}^{\top}C\left(s,\tau\right)\lambda-\lambda^{\top}\left(z-A\left(s,\tau\right)x-b\left(s,\tau\right)\right).
\label{eq:LQT_V_form2}
\end{align}
Here $\eta\left(s,\tau\right)$ and $b\left(s,\tau\right)$ are $n_x$-dimensional vectors, $A\left(s,\tau\right)$ is an  $n_x\times n_x$ matrix, $J\left(s,\tau\right)$ and $C\left(s,\tau\right)$ are symmetric, positive semi-definite $n_x\times n_x$ matrices. Function $\zeta (s,\tau)$ does not depend on the states, and does not need to be computed. The form of the conditional value function in \eqref{eq:LQT_V_form1} is equivalent to the discrete case \cite{Sarkka:2022} and corresponds to the dual problem of a quadratic program with affine equality constraints \cite{Boyd_book04}. 

It is direct to check that for a zero-length interval, the conditional value function \eqref{eq:V_zero_length} is recovered by setting $A\left(\tau,\tau\right)=I$, $b\left(\tau,\tau\right)=0$, $C\left(\tau,\tau\right)=0$, $\eta\left(\tau,\tau\right)=0$, and $J\left(\tau,\tau\right)=0$. In addition, if $C\left(s,\tau\right)$ is invertible, one can solve \eqref{eq:LQT_V_form1} in closed form and obtain a quadratic expression for $V(x,s;z,\tau)$ (cf.~Eq.~(34) in \cite{Sarkka:2022}):
\begin{align}
V(x,s;z,\tau)
&=\zeta(s,\tau)+\frac{1}{2}x^{\top}J\left(s,\tau\right)x-x^{\top}\eta\left(s,\tau\right) \nonumber \\
&\quad +\frac{1}{2}\left(z-A\left(s,\tau\right)x-b\left(s,\tau\right)\right)^\top \, C\left(s,\tau\right)^{-1} \nonumber \\
&\qquad \times \left(z-A\left(s,\tau\right)x-b\left(s,\tau\right)\right). \label{eq:LQT_V_form2b}
\end{align}
That is, if $C\left(s,\tau\right)$ is invertible, the conditional value function has the form \eqref{eq:LQT_V_form2b}, which is quadratic in $x$ and $z$. However, in general, $C\left(s,\tau\right)$ cannot be inverted, so it is suitable to use the dual function representation in \eqref{eq:LQT_V_form1}, which is more general. Matrix $C\left(s,\tau\right)$ fails to be invertible when some values $z(\tau)$ are not reachable from $x(s)$ implying that $V(x,s;z,\tau) = \infty$ which in LQT case translates to non-invertibility of $C\left(s,\tau\right)$. Because we are not assuming complete controllability of the system, this kind of state pairs can well exist even when the value function $V(x,t) < \infty$ for all $x$ and $t$.

The associative operator $\otimes$ in Theorem \ref{the:HJB-comb} for conditional value functions of the form \eqref{eq:LQT_V_form1} results in the following lemma which is analogous to \cite[Lemma 10]{Sarkka:2022}.
\begin{lem} \label{lem:LQT_combination}
	Given $V(x,s;z,\tau)$ and $V(z,\tau;y,t)$ of the form \eqref{eq:LQT_V_form1} with $s<\tau<t$, the conditional value function 
	\begin{equation}
	\begin{split}
	V(x, s ; y, t) &=V(x, s ; \cdot, \tau) \otimes V(\cdot, \tau ; y, t),
	\end{split}
	\end{equation}
	which is obtained using Theorem \ref{the:HJB-comb}, is of the form \eqref{eq:LQT_V_form1} with parameters
	\begin{equation}
	\begin{split}
	A\left(s,t\right) &= A\left(\tau,t\right) (I + C\left(s,\tau\right) J\left(\tau,t\right))^{-1} A\left(s,\tau\right), \\
	b\left(s,t\right) &= A\left(\tau,t\right) (I + C\left(s,\tau\right) J\left(\tau,t\right))^{-1} \\
	&	\quad \times(b\left(s,\tau\right) + C\left(s,\tau\right) \eta\left(\tau,t\right)) + b\left(\tau,t\right), \\
	C\left(s,t\right) &= A\left(\tau,t\right) (I + C\left(s,\tau\right) J\left(\tau,t\right))^{-1} C\left(s,\tau\right) A\left(\tau,t\right)^\top  \\
	&\quad + C\left(\tau,t\right), \\
	\eta\left(s,t\right) &= A\left(s,\tau\right)^\top (I + J\left(\tau,t\right) C\left(s,\tau\right))^{-1} \\
	&\quad	\times (\eta\left(\tau,t\right) - J\left(\tau,t\right) b\left(s,\tau\right)) + \eta\left(s,\tau\right), \\
	J\left(s,t\right) &= A\left(s,\tau\right)^\top (I + J\left(\tau,t\right) C\left(s,\tau\right))^{-1} J\left(\tau,t\right) A\left(s,\tau\right) \\
	&\quad + J\left(s,\tau\right),
	\end{split}
	\label{eq:lqt_comb}
	\end{equation}
	where $I$ is an identity matrix of size $n_x$.
\end{lem}
Because $C(s;\tau)$ and $J(\tau;t)$ are symmetric and positive semi-definite, it follows from Lemma~\ref{lem:ijc} in Appendix~\ref{sec:Invertible_matrices} that matrices $I+C(s;\tau)J(\tau;t)$ and $I+J(\tau;t)C(s;\tau)$ are always invertible. Lemma \ref{lem:LQT_combination} can be used to associatively combine the conditional value function parameters $A$, $b$, $C$, $\eta$, and $J$ once we know them for some non-vanishing time intervals. To obtain them for non-vanishing time intervals, we next solve the conditional HJB equations in Section \ref{subsec:Cond_back_forw_HJB} for the LQT case.

\subsection{Backward conditional HJB equation for LQT}\label{sub:conditional_backward_LQT}
This section explains how to solve the backward conditional HJB equation in Section \ref{subsec:Cond_back_forw_HJB} for the LQT problem. The resulting expressions will be useful to calculate the associative elements, which will be explained in Section \ref{sub:parallel_LQT}.
\begin{thm}
	\label{thm:Backward_HJB_LQT}
	For the LQT problem \eqref{eq:LQT_Control_problem} and a conditional value function $V(x,s;z,\tau)$ of the form \eqref{eq:LQT_V_form1}, the derivative  $\frac{\partial V(x,s;z,\tau)}{\partial s}$, given by the backward conditional HJB equation \eqref{eq:backward_conditional_HJB}, can be written in terms of the derivatives of the conditional value function parameters as
	\begin{equation}
	\begin{split}
\frac{\partial A\left(s,\tau\right)}{\partial s}
&= A\left(s,\tau\right)L\left(s\right)U\left(s\right)^{-1}L\left(s\right)^{\top}J\left(s,\tau\right) \\
&\quad
-A\left(s,\tau\right)F\left(s\right), \\
\frac{\partial b\left(s,\tau\right)}{\partial s}
&=-A\left(s,\tau\right)L\left(s\right)U\left(s\right)^{-1}L\left(s\right)^{\top}\eta\left(s,\tau\right) \\
&\quad-A\left(s,\tau\right)c\left(s\right), \\
\frac{\partial C\left(s,\tau\right)}{\partial s}
&=-A\left(s,\tau\right)L\left(s\right)U\left(s\right)^{-1}L\left(s\right)^{\top}A\left(s,\tau\right)^{\top}, \\ 
\frac{\partial\eta\left(s,\tau\right)}{\partial s}
&=-H\left(s\right)^{\top}X\left(s\right)r\left(s\right)-F\left(s\right)^{\top}\eta\left(s,\tau\right) \\
&\quad+J\left(s,\tau\right)c\left(s\right) \\
&\quad+J\left(s,\tau\right)L\left(s\right)U\left(s\right)^{-1}L\left(s\right)^{\top}\eta\left(s,\tau\right), \\
\frac{\partial J\left(s,\tau\right)}{\partial s}
&=-H\left(s\right)^{\top}X\left(s\right)H\left(s\right)-J\left(s,\tau\right)F\left(s\right) \\
&\quad-F\left(s\right)^{\top}J\left(s,\tau\right) \\
&\quad+J\left(s,\tau\right)L\left(s\right)U\left(s\right)^{-1}L\left(s\right)^{\top}J\left(s,\tau\right). \\
	\end{split}
	\label{eq:Backward_HJB_LQT}
	\end{equation} 
\end{thm}
\vspace*{0.2cm}

The proof of Theorem \ref{thm:Backward_HJB_LQT} and the proof that the conditional value function is of the form \eqref{eq:LQT_V_form1} are given in Appendix \ref{sec:Proof_backward_LQT_append}. This theorem enables us to solve the backward conditional HJB equation \eqref{eq:backward_conditional_HJB} in closed form for the LQT case. The required boundary conditions are
\begin{equation}
\begin{split}
  A\left(t,t \right) = I, \quad
  b\left(t,t \right ) &=0, \quad
  C\left(t,t \right) = 0, \\
  \eta\left(t,t \right) &= 0, \quad
  J\left(t,t \right) = 0, \\
\end{split}
\label{eq:abcej_boundary}
\end{equation} 
for all $t \in [t_0,t_f]$. 
It is also worth noting that the equations for $\eta$ and $J$ have exactly the same form as equations \eqref{eq:LQT_differential_S_v} for $v(t)$ and $S(t)$. This is no coincidence, because we indeed have $\eta(s,t_f^+) = v(s)$ and $J(s,t_f^+) = S(s)$, but we will come back to this connection in Section~\ref{sub:parallel_LQT}. However, this also implies that the solutions to them exist under the same conditions as they exists for $v(t)$ and $S(t)$, that is, under the conditions listed under \eqref{eq:LQT_Control_problem}.

\subsection{Forward conditional HJB equation for LQT}\label{sub:conditional_forward_LQT}
This section explains how to solve the forward conditional HJB equation in Section \ref{subsec:Cond_back_forw_HJB} for the LQT problem. The resulting expressions will be useful to calculate the associative elements, which will be explained in Section \ref{sub:parallel_LQT}.
\begin{thm}
	\label{thm:Forward_HJB_LQT}
	For the LQT problem \eqref{eq:LQT_Control_problem} and a conditional value function $V(x,s;z,\tau)$ of the form \eqref{eq:LQT_V_form1}, the derivative  $\frac{\partial V(x,s;z,\tau)}{\partial\tau}$, given by the forward conditional HJB equation \eqref{eq:forward_conditional_HJB}, can be written in terms of the derivatives of the conditional value function parameters such that
	\begin{equation}
	\begin{split}
\frac{\partial A\left(s,\tau\right)}{\partial\tau}
&=-C\left(s,\tau\right)H\left(\tau\right)^{\top}X\left(\tau\right)H\left(\tau\right)A\left(s,\tau\right) \\
&\quad+F\left(\tau\right)A\left(s,\tau\right), \\
\frac{\partial b\left(s,\tau\right)}{\partial\tau}
&=C\left(s,\tau\right)H\left(\tau\right)^{\top}X\left(\tau\right)r\left(\tau\right)+F\left(\tau\right)b\left(s,\tau\right) \\
&\quad-C\left(s,\tau\right)H\left(\tau\right)^{\top}X\left(\tau\right)H\left(\tau\right)b\left(s,\tau\right)+c\left(\tau\right), \\
\frac{\partial C\left(s,\tau\right)}{\partial\tau}
&=-C\left(s,\tau\right)H\left(\tau\right)^{\top}X\left(\tau\right)H\left(\tau\right)C\left(s,\tau\right) \\
&\quad+L\left(\tau\right)U\left(\tau\right)^{-1}L\left(\tau\right)^{\top}+F\left(\tau\right)C\left(s,\tau\right) \\
&\quad+C\left(s,\tau\right)F\left(\tau\right)^{\top}, \\
\frac{\partial\eta\left(s,\tau\right)}{\partial\tau}
&=A\left(s,\tau\right)^{\top}H\left(\tau\right)^{\top}X\left(\tau\right)r\left(\tau\right) \\
&\quad-A\left(s,\tau\right)^{\top}H\left(\tau\right)^{\top}X\left(\tau\right)H\left(\tau\right)b\left(s,\tau\right), \\
\frac{\partial J\left(s,\tau\right)}{\partial\tau}
&=A\left(s,\tau\right)^{\top}H\left(\tau\right)^{\top}X\left(\tau\right)H\left(\tau\right)A\left(s,\tau\right).
	\end{split}
	\label{eq:Forward_HJB_LQT}
	\end{equation}
\end{thm}
\vspace*{0.2cm}
Theorem \ref{thm:Forward_HJB_LQT} and that the conditional value function is of the form \eqref{eq:LQT_V_form1} are proved in Appendix \ref{sec:Proof_forward_LQT_append}. The boundary conditions for the parameters are given by \eqref{eq:abcej_boundary} for $t \in [t_0,t_f]$. Because these solutions are the same as those of \eqref{eq:Backward_HJB_LQT}, the existence conditions for the solutions are also the same.

\subsection{Parallel solution of optimal LQT control}\label{sub:parallel_LQT}
To obtain the parallel solution of optimal LQT control using the parallel scan algorithm, we first need to compute the elements (conditional value functions), see Theorem \ref{the:parallel_HJB}. The elements $a_{j}$ for $j\in\{0,\ldots,T-1\}$ in \eqref{eq:a_j} where $V(\cdot,t_{j};\cdot,t_{j+1})$ is of the form \eqref{eq:LQT_V_form1}. can be computed by solving the conditional backward or forward HJB equations, given by Theorems \ref{thm:Backward_HJB_LQT} and \ref{thm:Forward_HJB_LQT}, with the boundary conditions given in \eqref{eq:abcej_boundary}.
The element 
\begin{equation}
a_{T}=V(x,t_{f})=V(x,t_{f};0,t_{f}^+),
\end{equation}
where $V(x,t_{f};0,t_{f}^+)$ is of the form \eqref{eq:LQT_V_form1} with parameters
\begin{equation}\label{eq:LQT_final_initialisation}
\begin{split}
A\left(t_{f},t_{f}^{+}\right) = 0, \quad
b\left(t_{f},t_{f}^{+}\right) &= 0, \quad
C\left(t_{f},t_{f}^{+}\right) = 0, \\
J\left(t_{f},t_{f}^{+}\right) &= H_f^\top \, X_f \, H_f, \\
\eta\left(t_{f},t_{f}^{+}\right) &= H_f^\top \,  X_f \, r_f.
\end{split}
\end{equation}

We should note that if we impose these terminal conditions to \eqref{eq:Backward_HJB_LQT}, then we will get the solutions $A\left(t,t_{f}^+\right)=0$, $b\left(t,t_{f}^+\right)=0$ and $C\left(t,t_{f}^+\right)=0$ for $t \in [t_0,t_f]$. Furthermore, the equations for $\eta$ and $J$ then become
\begin{equation}
\begin{split}
\frac{\partial\eta\left(t, t_{f}^+\right)}{\partial t}
&=-H\left(t\right)^{\top}X\left(t\right)r\left(t\right)-F\left(t\right)^{\top}\eta\left(t, t_{f}^+\right) \\
&\quad
+J\left(t, t_{f}^+\right)c\left(t\right) \\
&\quad+J\left(t, t_{f}^+\right)L\left(t\right)U\left(t\right)^{-1}L\left(t\right)^{\top}\eta\left(t, t_{f}^+\right), \\
\frac{\partial J\left(t, t_{f}^+\right)}{\partial t}
&=-H\left(t\right)^{\top}X\left(t\right)H\left(t\right)-J\left(t, t_{f}^+\right)F\left(t\right) \\
&\quad-F\left(t\right)^{\top}J\left(t,t_{f}^+\right) \\
&\quad+J\left(t, t_{f}^+\right)L\left(t\right)U\left(t\right)^{-1}L\left(t\right)^{\top}J\left(t, t_{f}^+\right),
\end{split}
\end{equation}
which together with \eqref{eq:LQT_final_initialisation}, imply that, as expected, we recover $v(t)$ and $S(t)$ in \eqref{eq:LQT_differential_S_v} by
\begin{equation}
\begin{split}
	\eta\left(t,t_{f}^+\right) = v(t), \quad J\left(t,t_{f}^+\right) = S(t).
\end{split}
\label{eq:eJ_vS}
\end{equation}

After we have computed all elements $a_{j}$, we can use the parallel scans in Theorem \ref{the:parallel_HJB} to obtain the parameters of the conditional value functions $V(\cdot,t_0;\cdot,t_k)$ and value functions $V(\cdot,t_{k})$ in parallel. To obtain the conditional value function, or the value function, at any given time within each time sub-interval, we can solve the ODEs in either of Theorems \ref{thm:Backward_HJB_LQT} or \ref{thm:Forward_HJB_LQT}. This can be done in parallel for each sub-interval. In addition, within each sub-interval, it is possible to parallelise operations by splitting the sub-interval into two equal sub-intervals, using the forward HJB equation for the first sub-interval, and the backward HJB equation for the second sub-interval.

Once we have computed $V(x,t)$ for all times, we can compute the optimal control using \eqref{eq:LQT_optimal_control}, which can be done in parallel for all $t$. 

\subsection{Optimal trajectory recovery}\label{sub:optimal_trajectory_LQT}
We proceed to explain Methods 1 and 2 presented in Section \ref{subsec:Optimal_trajectory} for optimal trajectory recovery, for the LQT case.

\subsubsection{Method 1}
In the LQT case, the optimal trajectory can be computed using \eqref{eq:optimal_trajectory_LQT}. Then, for each sub-interval, we can obtain $\Psi(t,t_{j-1})$ and $\alpha(t,t_{j-1})$ by solving the ODEs \eqref{eq:ODE_transition_matrix} and \eqref{eq:ODE_phi} in parallel. If we have $\Psi(\tau,s)$, $\alpha(\tau,s)$ in an interval $\left[s,\tau\right]$, and $\Psi(t,\tau)$, $\alpha(t,\tau)$ in an interval $\left[\tau,t\right]$, we can obtain $\Psi(t,s)$ and $\alpha(t,s)$ using the combination rule
\begin{align}
\Psi(t,s) &= \Psi(t,\tau) \, \Psi(\tau,s), \label{eq:Psi_comb} \\
\alpha(t,s) &= \Psi(t,\tau)\, \alpha(t,\tau) + \alpha(\tau,s).
\label{eq:alpha_comb}
\end{align}
This combination rule is associative \cite[Sec. IV.C]{Sarkka_book19}. Then, we can compute $\Psi(t_j,t_{j-1})$ and $\alpha(t_j,t_{j-1})$ for all $j$ and combine them using parallel scans to obtain $\Psi(t_j,t_{0})$ and $\alpha(t_j,t_{0})$ for all $j$.

For $t\in (t_{j-1},t_{j})$, the optimal trajectory is given by
\begin{align}
x^{*}(t) &= \Psi(t,t_{j-1}) \, \Psi(t_{j-1},t_{0}) \, x\left(t_{0}\right)\nonumber \\
& \quad +\Psi(t,t_{j-1}) \, \alpha(t_{j-1},t_{0}) + \alpha(t,t_{j-1}),
\end{align}
which can be computed in parallel for each sub-interval.

\subsubsection{Method 2}
If $V(x,t)$ is of the form \eqref{eq:V_quadratic_LQT} and $V(x(t_{0}),t_{0},x,t)$ is of the form \eqref{eq:LQT_V_form1}, the solution to \eqref{eq:Opt_trajectory_method2} is \cite[Lemma 15]{Sarkka:2022}
\begin{align}
x^*(t) & = (I + C(t_0,t) \, S(t))^{-1} \nonumber \\
&\quad \times (A(t_0,t) \, x(t_0) + b(t_0,t) + C(t_0,t) \, v(t)).
\end{align}
This operation can be done in parallel for all the discretised time steps. Also recall from \eqref{eq:eJ_vS} that we actually have $v(t) = \eta\left(t,t_{f}^+\right)$ and $S(t) = J\left(t,t_{f}^+\right)$.

\section{Computing conditional value functions using direct optimisation} \label{sec:direct_opt}

The direct use of finite difference schemes \cite{Crandall84} for solving the conditional value functions in \eqref{eq:backward_conditional_HJB} and \eqref{eq:forward_conditional_HJB} with boundary condition \eqref{eq:V_zero_length} can be numerically challenging. Instead, we can recall that a conditional value function solves an optimal control problem with a terminal constraint on the state (see Definition \ref{Defn:conditional_value_function}). Then, at least in principle, we can use any optimal continuous-time control solver that admits terminal constraints to compute conditional value functions. For example, one can use solutions based on calculus of variations to obtain the equations for optimal controls and trajectories \cite[Chap. 3]{Lewis_book12}. In some cases, closed-form results can be obtained \cite[Chap. 3]{Lewis_book12}, and otherwise, we can use numerical methods such as direct and indirect shooting methods to obtain approximate solutions \cite{Bock84,Betts98}.

In the experimental results in the next section, we make use of direct shooting methods \cite{Bock84,Betts98}. The idea is to partition the state space to a grid of points and calculate the conditional value functions for each possible start and end point on the grid. Let us consider the calculation of the element $a_{j}=V(\cdot,t_{j-1};\cdot,t_{j})$ for $j\in\{1,\cdots,T-1\}$ for given start and end points. To apply a direct shooting method, we divide the time interval into $n$ smaller intervals $t_{j}^{0}<t_{j}^{1}<...<t_{j}^{n}$ where $t_{j}^{0}=t_{j-1}$ and $t_{j}^{n}=t_{j}.$ At each time interval we form a parametric approximation of the control
\begin{equation}
u(t)\approx\phi_{j}^{m}(t,q_{j}^{m}),\quad t\in[t_{j}^{m},t_{j}^{m+1}),
\end{equation}
where $\phi_{j}^{m}(\cdot,\cdot)$ is a parametric function with parameters $q_{j}^{m}$. We then approximate the solution from $x(t_{j}^{m})=s_{j}^{m}$
to $x(t_{j}^{m+1})$, which is denoted as $x(t;s_{j}^{m},q_{j}^{m})$, by numerically solving
\begin{equation}
\begin{split}\frac{dx(t)}{dt} & =f(x(t),\phi_{j}^{m}(t,q_{j}^{m}),t),\quad x(t_{j}^{m})=s_{j}^{m}\end{split},
\end{equation}
which is a function of $q_{j}^{m}$ and $s_{j}^{m}$. The cost in the interval $[t_{j-1},t_{j}]$ can then be approximated as 
\begin{equation}
\begin{split}C_{j}(q_{j},s_{j})\approx\sum_{m=0}^{n-1}\int_{t_{j}^{m}}^{t_{j}^{m+1}}\ell(x(t),\phi_{j}^{m}(t,q_{j}^{m}),t),t)\,dt\end{split}
\end{equation}
which is to be minimised with respect to $q_{j}=(q_{j}^{0},...,q_{j}^{n-1})$ and $s_{j}=(s_{j}^{0},\ldots,s_{j}^{n-1})$ with the constraints
\begin{equation}
\begin{split}x(t_{j}^{m+1};s_{j}^{m},q_{j}^{m})-s_{j}^{m+1} & =0.\end{split}
\end{equation}
We can also ensure the continuity of the control with the constraints
\begin{equation}
\phi_{j}^{m}(t_{j}^{m+1},q_{j}^{m})-\phi_{j}^{m+1}(t_{j}^{m+1},q_{j}^{m+1})=0.
\end{equation}
Therefore, together with the terminal constraint $x(t_{j})=z$, we obtain an optimisation problem with equality constraints. This type of problem can be solved, for example, by using sequential quadratic programming \cite[Chap. 18]{Nocedal+Wright:2006}. The required derivatives can be easily computed with automatic differentiation.

Once we have obtained the conditional value functions, we can combine them via Theorem~\ref{the:HJB-comb} by directly minimising over the grid of points. In practice, the direct use of Theorem~\ref{the:HJB-comb} on the grid is numerically inaccurate. Fortunately, we can refine the numerical accuracy of this step by using interpolation, for example, by forming a local quadratic approximation of the conditional value function, and then performing the minimisation. The interpolation we use in the nonlinear control experiment in Section \ref{subsec:nonlinear_control_experiment} is explained in Appendix \ref{sec:Interpolation_append}. This interpolation though assumes that the conditional value functions are twice continuously differentiable in the state variable.

\section{Numerical results} \label{sec:results}

In this section, we present proof-of-concept results of the proposed algorithms both on a multi-core CPU and a GPU. The purpose of these experiments is to show that the methods provide real computational benefit at least in a simple setting although more comprehensive evaluation is needed in the future. For experimental evaluations we used TensorFlow 2.10 \cite{Abadi_et_al:2015} parallel computing framework. The CPU experiments were run on AMD EPYC 7643 48-core processor with 512GB of memory. The GPU experiments were run using NVIDIA A100-SXM GPU with 80GB of memory.

\subsection{Continuous-time LQT} \label{sec:lqt_experiment}

The associative rules for combining general value functions for LQT problems allow for implementing parallel solver for the LQT by using the associative scan algorithm. The aim of this experiment is to empirically test the performance of this kind of methods on a multi-core CPU and a GPU. The model we consider is a continuous-time version of the LQT model considered in \cite{Sarkka:2022} which has the form
\begin{equation}
\begin{split}
  f(x,u) &= \underbrace{\begin{bmatrix}
    0 & 0 & 1 & 0 \\
    0 & 0 & 0 & 1 \\
    0 & 0 & 0 & 0 \\
    0 & 0 & 0 & 0
  \end{bmatrix}}_{F} \, x
  + \underbrace{\begin{bmatrix}
    0 & 0  \\
    0 & 0 \\
    1 & 0 \\
    0 & 1 
  \end{bmatrix}}_{L} \, u, 
\end{split}
\end{equation}
where the state contains the position and velocity $x = \begin{bmatrix} p_x & p_y & v_x & v_y \end{bmatrix}^\top$ and the control signals are the accelerations $u = \begin{bmatrix} a_x & a_y \end{bmatrix}^\top$. The cost function is given by
\begin{equation}
\begin{split}
\ell(x,u,t) &= \frac{1}{2} (r(t) - H \, x)^\top X (r(t) - H \, x) + \frac{1}{2} u^\top U \, u, \\
\phi(x) &= \frac{1}{2} (H_f \, x - r_f)^\top \, X_f \, (H_f \, x - r_f),
\end{split}
\end{equation}
where
\begin{equation}
  H = \begin{bmatrix}
    1 & 0 & 0 & 0 \\
    0 & 1 & 0 & 0
  \end{bmatrix}, \quad X = I_{2 \times 2}, \quad 
  U = 10^{-1} \, I_{2 \times 2},
\end{equation}
and $H_f = X_f = I_{4 \times 4}$. The length of the time interval is $t_f = 50$.

\begin{figure}[tb]
\centering
\includegraphics{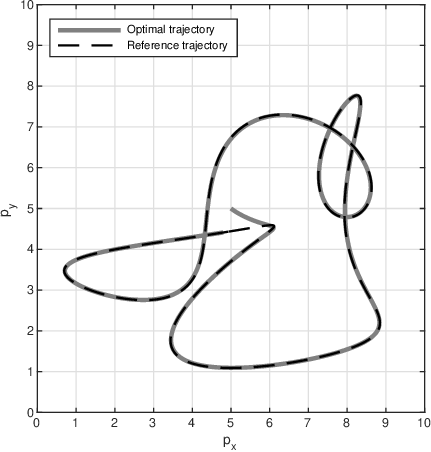}
\caption{The reference trajectory of the parallel LQT experiment along with the optimal trajectory. The diversion of the optimal path from the reference trajectory in the middle is due to the starting point being slightly off the reference trajectory.}
\label{fig:track_trajectory}
\end{figure}

The function $r(t)$ is a function which returns the 2d trajectory depicted in Fig.~\ref{fig:track_trajectory}. The trajectory was generated by fitting a Fourier series to a sampled trajectory from the discrete-time model considered in \cite{Sarkka:2022}. The initial state is $x(0) =  \begin{bmatrix} 5 & 5 & 0 & 0 \end{bmatrix}^\top$ which is located slightly off from the starting of the reference trajectory. Thus the aim of the control problem is to follow this reference trajectory as accurately as possible while consuming a minimal amount of energy. The optimal trajectory for this model is also shown in the same figure. 

We implemented sequential and parallel solvers for the backward pass which computes the value function parameters $S(t)$ and $v(t)$, and the forward pass which computes the optimal control and trajectory $x^*(t)$ and $u^*(t)$. For this purpose, we split the time span of the problem into $T$ blocks with $n$ substeps in each of the blocks. We used $n = 10$ in the experiment and varied the number of blocks $T$. In the sequential solution method, the value function parameters were solved by applying the 4th order Runge--Kutta method \cite{Sarkka_book19} to the backward equations \eqref{eq:LQT_differential_S_v} on the grid of $n \, T$ points on the interval $[0,t_f]$. The optimal trajectory was then computed by solving for the control law $u(x,t)$ and integrating the forward equations \eqref{eq:optimal_trajectory} by using the 4th order Runge--Kutta method. The complexity of this sequential solution is $O(T)$.

In the parallel method, we first solved the conditional value functions for each block in parallel by using $n = 10$ steps of 4th order Runge--Kutta integration for the backward equations in Theorem~\ref{thm:Backward_HJB_LQT} (we could have equivalently used Theorem~\ref{thm:Forward_HJB_LQT}). The blocks were combined by applying the backward parallel associative scan algorithm to the blocks with the associative operator defined in Lemma~\ref{lem:LQT_combination}. Finally, the intermediate values within each block were computed by solving (in parallel) the backward equations \eqref{eq:LQT_differential_S_v} using $n$ steps of 4th order Runge--Kutta. The (span) complexity of this parallel solution is $O(\log T)$.

In the parallel method, the optimal trajectory was solved using the two different methods described in Sec.~\ref{sub:optimal_trajectory_LQT}. In Method 1, we first solved $\Psi$ and $\alpha$ for each block from equations \eqref{eq:ODE_transition_matrix} and \eqref{eq:ODE_transition_matrix2} in parallel and then applied parallel associative scan to operator defined by \eqref{eq:Psi_comb} and \eqref{eq:alpha_comb}. The final result was computed by solving \eqref{eq:optimal_trajectory_ODE_LQT} for each block, again in parallel. In Method 2, the block conditional value functions were computed similarly to the parallel backward pass although this time using the equations from Theorem~\ref{thm:Forward_HJB_LQT} (Theorem~\ref{thm:Backward_HJB_LQT} could have been used as well). The blocks were combined using parallel associative scan with operator defined in Lemma~\ref{lem:LQT_combination}. The intermediate values of $A$, $b$, and $C$ inside the blocks were solved from their equations given in Theorem~\ref{thm:Forward_HJB_LQT} after which the optimal trajectory was solved using \eqref{eq:Opt_trajectory_method2}. All the differential equation solutions were computed by using the 4th order Runge--Kutta method. The (span) complexity of this procedure is $O(\log T)$.

\begin{figure}[tb]
\centering
\includegraphics{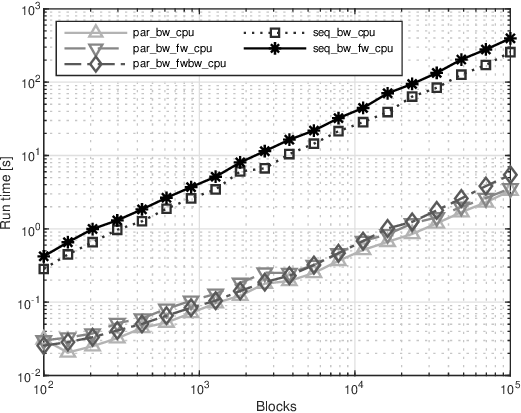}
\caption{Results of parallel continuous-time LQT experiment ran on CPU. The parallel algorithms (prefix 'par') can be seen to be consistently faster than the sequential ones (prefix 'seq').}
\label{fig:clqt_cpu}
\end{figure}

The results of the experiment run on the CPU are shown in Fig.~\ref{fig:clqt_cpu}. The parallel methods consisting of backward pass only is labeled as {\em par\_bw\_cpu}, the method consisting of backward pass followed by Method 1 based forward pass is labeled as {\em par\_bw\_fw\_cpu}, and the method consisting of backward pass followed by optimal trajectory computation method with Method 2 is labeled as {\em par\_bw\_fwbw\_cpu}. The sequential versions of the backward pass only and the combined backward and forward passes are labeled as {\em seq\_bw\_cpu} and {\em par\_bw\_fw\_cpu}, respectively. It can be seen that the parallel methods are consistently faster than the sequential methods with a quite constant multiplicative factor of close to two decades. 

\begin{figure}[tb]
\centering
\includegraphics{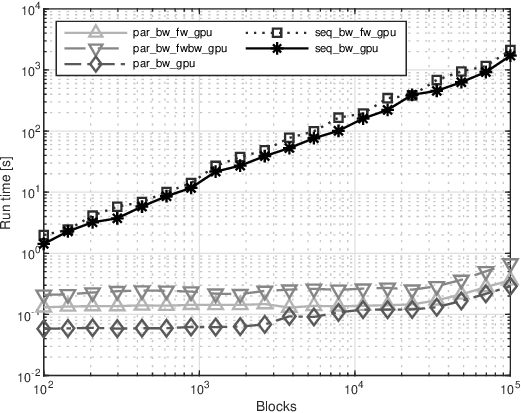}
\caption{Results of parallel continuous-time LQT experiment ran on GPU. It can be seen that the parallel algorithms (prefix 'par') are orders of magnitude faster than the sequential algorithms (prefix 'seq').}
\label{fig:clqt_gpu}
\end{figure}

Fig.~\ref{fig:clqt_gpu} shows the corresponding results for the GPU. The different methods are labeled using the same logic as in the CPU experiment above. It can be seen that until up to $10^4$ blocks, the run time of the parallel methods is almost constant as compared to the linear growth of the cost of the sequential methods. Only closer to $10^5$ blocks there are signs of an emerging linear growth while still the difference to the sequential method is huge. Although it would be interesting to investigate the behaviour with larger block counts, unfortunately, the amount of available GPU memory did not allow for that.

\begin{figure}[tb]
\centering
\includegraphics{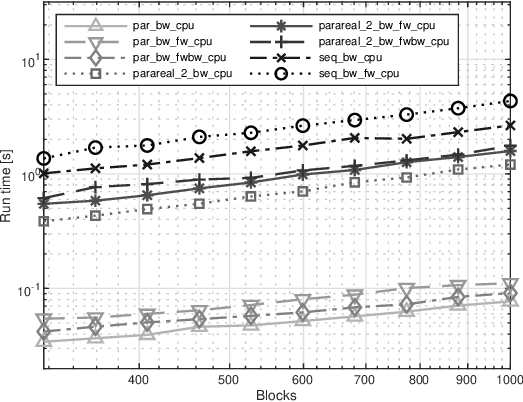}
\caption{Results of comparison of parareal and associative scan methods for continuous-time LQT run on CPU. The parareal algorithm can be seen to be faster than the sequential algorithms, but slower than the parallel algorithms.}
\label{fig:parareal_clqt_cpu}
\end{figure}

\begin{figure}[tb]
\centering
\includegraphics{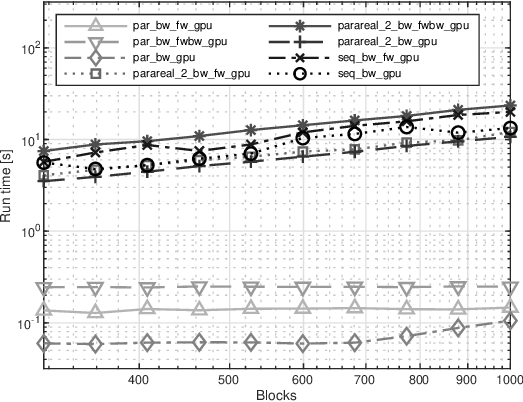}
\caption{Results of comparison of parareal and associative scan methods for continuous-time LQT run on GPU. The parareal algorithm seems to have roughly the same speed as the sequential algorithms while the proposed parallel algorithms are an order of magnitude faster.}
\label{fig:parareal_clqt_gpu}
\end{figure}

For comparison, we also implemented the parareal method \cite{Lions:2001,gander2007analysis} for solving the backward and forward equations \eqref{eq:LQT_differential_S_v}  and \eqref{eq:optimal_trajectory}. Due to challenging implementation and limitations in memory limitations, we could only run parareal up to 1\,000 blocks, but the trend is already visible in the results. The parareal is also numerically stable only for short enough steps in the coarse integration method and hence we had to start the number of blocks from around 300 instead of 100 which we used above. We used 2 iterations in parareal which seemed to be enough to give a quite accurate result already. The results run on the CPU are shown in Fig.~\ref{fig:parareal_clqt_cpu} and the results run on GPU are shown in Fig.~\ref{fig:parareal_clqt_gpu}. On the CPU, parareal provides and advantage over the sequential methods. However, the proposed parallel methods are faster than it by a clear margin. On the GPU, parareal has roughly the same performance as the sequential methods.

\subsection{Non-linear control problem}\label{subsec:nonlinear_control_experiment}

In this experiment the aim is to illustrate the methodology in a non-linear model which we can approximate on a finite grid using finite differences and where the conditional value functions can be evaluated by a classical optimal control solution method, a direct shooting method in this case. For that purpose, we consider the scalar model
\begin{equation}
\begin{split}
  \frac{dx}{dt} = f(x,u) = 1 - \beta \, x^2 + u, 
\end{split}
\label{eq:nonlin_ex_f}
\end{equation}
where $\beta = 1/10$, which is a simple model of the velocity of a falling body under air resistance which we can control with force $u$. The cost function is
\begin{equation}
\begin{split}
  \ell(x,u) &= \frac{1}{2} \,  x^2 + \frac{1}{2} \, u^2, \\
  \phi(x) &= 2 \, x^2,
\end{split}
\label{eq:nonlin_ex_ell_phi}
\end{equation}
that is, the aim is to (approximately) stop the fall while limiting the control energy. In this experiment we fix the block length to $1/10$ time units and vary the number of blocks $T$ to avoid numerical issues from the varying block length. Each block is further divided to $n=10$ steps both in the sequential and parallel methods to aid numerical accuracy although the final aim is to solve the HJB equation \eqref{eq:HJB_equation}, with $f$, $\ell$, and $\phi$ as defined in \eqref{eq:nonlin_ex_f} and \eqref{eq:nonlin_ex_ell_phi}, on the time resolution of $1/10$. For the sequential solution we use the upwind discretisation method considered in \cite{Crandall84} where the time discretisation the block length divided by the number of steps $n=10$ while the results are only stored at the block resolution of $1/10$ time units.

In the parallel solution we compute the conditional value functions for the blocks using a direct shooting method (see Sec.~\ref{sec:direct_opt}) over the $n=10$ steps. In the shooting method, we use a piecewise linear approximation for the control and the integration of the state and the cost are done using the 4th order Runge--Kutta method. The cost function is optimized by computing the first and section derivatives of the incurred cost by using automatic differentation in TensorFlow and by solving the resulting constrained optimisation problem using local sequential quadratic programming (SQP) \cite[Chap. 18]{Nocedal+Wright:2006} (see also Sec.~\ref{sec:direct_opt}). Because the model is time-invariant, the conditional value function is the same for each block and hence only needs to be solved only for one block. This computation would though be fully parallelisable over the blocks. The block solutions were combined by running parallel associative scan which uses the interpolated combination operation described in Appendix~\ref{sec:Interpolation_append}. Hence, we assume that the conditional value functions are twice differentiable or can at least be accurately approximated as such.

\begin{figure}[tb]
\centering
\includegraphics{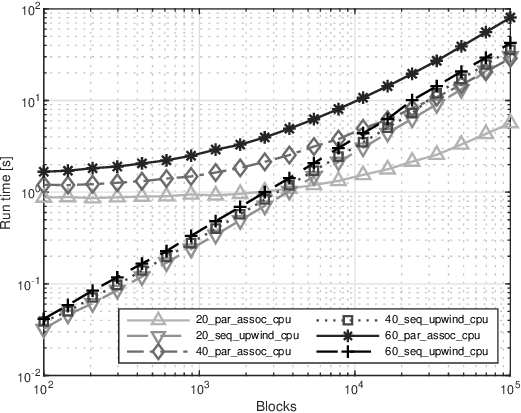}
\caption{Results of the experiment on non-linear sequential (upwind) and parallel (associative) control methods using grid sizes of 20, 40, and 60 run on the CPU.}
\label{fig:nonlin_cpu}
\end{figure}

Fig.~\ref{fig:nonlin_cpu} shows the result of running the TensorFlow implementation sequential (upwind) solution and the parallel solution described above on the CPU. The TensorFlow version and the CPU were the same as in experiment described in Section~\ref{sec:lqt_experiment}. The experiment was performed using grids of sizes 20, 40, and 60. It can be seen that the initial SQP method causes a large overhead to the parallel results and that the performance of the parallel method is heavily dependent on the grid size. With small grid sizes the parallel method is faster than the sequential method when the number of blocks increases, but with larger grid sizes, the sequential method appears to be faster also in this regime.

\begin{figure}[tb]
\centering
\includegraphics{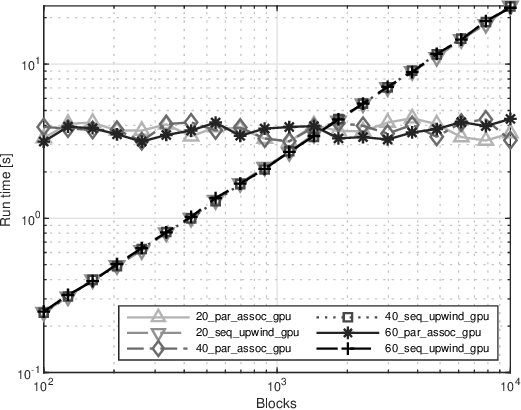}
\caption{Results of the experiment on non-linear sequential (upwind) and parallel (associative) control methods using grid sizes of 20, 40, and 60 run on the GPU.}
\label{fig:nonlin_gpu}
\end{figure}

The results of the sequential (upwind) solution and the parallel solution when run on GPU (the same as described in Sec.~\ref{sec:lqt_experiment}) are shown in Fig.~\ref{fig:nonlin_gpu}. Due to limitations in the GPU memory we were only able run the experiment up to $10^4$ blocks, but the trend is already clearly visible. The sequential methods have a consistent linear growth in run time, while the parallel methods only show a slight increase with the block length. The large constant bias in the parallel run times can be explained by the time taken by the SQP-based solving of the conditional value function. 

\section{Conclusion} \label{sec:conclusion}

In this paper, we have presented methods for temporal parallelisation of continuous-time optimal control problems. To do so, we first introduced the concept of conditional value functions, and obtained the backward and forward conditional HJB equations to calculate them. Then, we showed how to combine two conditional value functions using an associative operator. The parallelisation in the method is based on dividing the time interval into sub-intervals, calculating the conditional value functions in each sub-interval, and combining the results via parallel scans. We have also developed methods for optimal trajectory recovery using parallel scans. We also derived the resulting equations for the LQT control problem, which admits closed-form formulas for both forward and backward conditional HJB equations. We also demonstrated the methodology with analytical examples and also experimentally evaluated the computational performance of the parallel methods on multi-core CPU and GPU.

The extension of the proposed framework to parallelising general continuous-time stochastic control problems \cite{Stengel_book94} is not straightforward. However, due to the certainty equivalence property of linear stochastic control problems \cite{Stengel_book94}, the temporal parallelisation of continuous-time stochastic LQT problems can be performed with the approach presented in this paper without any modifications.

The main limitation of the present approach is that it can be interpreted as trading lower time complexity for more extensive storage requirements. Instead of single time step value functions, we need to store the conditional value functions for the start and end point pairs of the sub-intervals. Furthermore, to really achieve the theoretical parallelisation, an extensively high number of computational cores is sometimes needed.

\appendices

\section{Derivation of forward conditional HJB equation}\label{sec:Proof_forward_HJB_append}

In this appendix we derive the forward conditional HJB equation \eqref{eq:forward_conditional_HJB}. We recall that
the conditional value function $V(y,s;z,\tau)$, see Definition \ref{Defn:conditional_value_function}, solves
\begin{align}
\frac{dx(t)}{dt} & =f(x(t),u(t),t), \label{eq:forward_cond_append1} \\
C[u(t)] & =\int_{s}^{\tau}\ell(x(t),u(t),t)\,dt, \label{eq:forward_cond_append2}
\end{align}
with $x(s)=y$ and $x(\tau)=z$. We can now rewrite the cost as
\begin{align}
C[u(t)] & =\int_{\tau}^{s} -\ell(x(t),u(t),t)\,dt, \label{eq:forward_cond_append2b}
\end{align}
and define $x(\tau)=z$ to be the start state and $x(s)=y$ the end state. Although formally, we could now write the HJB equation for this in variables $z$ and $\tau$, this would lead to an incorrect equation, because the HJB equation is not valid when $\tau > s$. To get the correct equation, we need to do a change of variables as follows.

To derive the forward conditional HJB equation, we time-reverse the problem \eqref{eq:forward_cond_append1}-\eqref{eq:forward_cond_append2} by setting $t'=-t$ such that
\begin{align}
\frac{dx(-t')}{d(-t')} & =f(x(-t'),u(-t'),-t'),\\
C[u(-t')] & =\int_{-s}^{-\tau}\ell(x(-t'),u(-t'),-t')\,d(-t'),\nonumber \\
& =\int_{-\tau}^{-s}\ell(x(-t'),u(-t'),-t')\,dt',
\end{align}
with $x(s)=y$ and $x(\tau)=z$.

We now define the time-reversed functions $\bar{x}(t')=x(-t')$, $\bar{u}(t')=u(-t')$, and
\begin{align}
\bar{f}(\bar{x}(t'),\bar{u}(t'),t')&=f(x(-t'),u(-t'),-t'),\\
\bar{\ell}(\bar{x}(t'),\bar{u}(t'),t')&=\ell(x(-t'),u(-t'),-t').
\end{align}
Substituting these expressions into \eqref{eq:forward_cond_append1}-\eqref{eq:forward_cond_append2} yields
\begin{align}
\frac{d\bar{x}(t')}{dt'} & =-\bar{f}(\bar{x}(t'),\bar{u}(t'),t'),\\
C[\bar{u}(t')] & =\int_{-\tau}^{-s}\bar{\ell}(\bar{x}(t'),\bar{u}(t'),t')\,dt',
\end{align}
where $\bar{x}(-s)=y$ and $\bar{x}(-\tau)=z$. 

The conditional value function for this time-reversed system is denoted by $\overline{V}(z,\tau';y,s')$  where $s'=-s$ and $\tau'=-\tau$. Provided that it is differentiable in $z$ and $\tau'$ at least in a viscosity sense, it satisfies the backward HJB equation
\begin{align}
& -\frac{\partial\overline{V}(z,\tau';y,s')}{\partial\tau'}\nonumber \\
& =\min_{\bar{u}}\left\{ \bar{\ell}(z,\bar{u},\tau')-\left(\frac{\partial\overline{V}(z,\tau';y,s')}{\partial z}\right)^{\top}\,\bar{f}(z,\bar{u},\tau')\right\} ,
\end{align}
Because the backward conditional value function must be the original value function backwards,
\begin{equation}
\overline{V}(z,\tau';y,s')=V(y,-s';z,-\tau')
\end{equation}
we can reverse back the time to the get the forward equation for the conditional value
function
\begin{align}
& \frac{\partial V(y,s;z,\tau)}{\partial\tau}\nonumber \\
& =\min_{\bar{u}}\left\{ \ell(z,u,\tau)-\left(\frac{\partial V(y,s;z,\tau)}{\partial z}\right)^{\top}\,f(z,u,\tau)\right\} .
\end{align}

This finishes the derivation of \eqref{eq:forward_conditional_HJB}.

\section{Invertible matrices in Lemma \ref{lem:LQT_combination}}\label{sec:Invertible_matrices}

This appendix contains a lemma that ensures that matrices $I+C(s;\tau)J(\tau;t)$ and $I+J(\tau;t)C(s;\tau)$ in Lemma \ref{lem:LQT_combination} are invertible. The lemma is the following.
\begin{lem} \label{lem:ijc}
	Let $J$ and $C$ be two real, symmetric, and positive semi-definite matrices of the same size.  Then matrices $I+JC$ and $I+CJ$ are invertible.
\end{lem}
\begin{proof}
Let $v$ be a non-zero eigenvector of $I+JC$ and $\lambda$ the corresponding eigenvalue. We then have
\begin{align}
(I+JC)v &= \lambda v\label{eq:IJC_append1},\\
JCv &= (\lambda-1)v\label{eq:IJC_append2}.
\end{align}
We multiply from the left by $(Cv)^{\mathsf{H}}$, where the superscript $\mathsf{H}$
denotes the conjugate transpose. Then, \eqref{eq:IJC_append2} becomes
\begin{equation}
v^{\mathsf{H}}CJCv=(\lambda-1)v^{\mathsf{H}}Cv.
\end{equation}

Due to the positive semi-definiteness of $J$ and $C$, we then have
$v^{\mathsf{H}}CJCv\geq0$ and $v^{\mathsf{H}}Cv\geq0$ . When $v^{\mathsf{H}}Cv\neq0$, this
implies that $\lambda$ is real and
\begin{equation}
\lambda=1+(v^{\mathsf{H}}CJCv)/(v^{\mathsf{H}}Cv)\geq1.
\end{equation}
When $v^{\mathsf{H}}Cv=0$, then we must have $Cv=0$. This is because we can
always write $C=U^{\top}U$ for some matrix $U$ which enables us to
write $v^{\mathsf{H}}Cv=(Uv)^{\mathsf{H}}(Uv)=||Uv||^{2}$, which is zero only if $Uv=0$,
hence $Cv=U^{\top}Uv=U^{\top}0=0$. Thus \eqref{eq:IJC_append1} reduces
to
\begin{align}
v &= \lambda v,
\end{align}
which implies that the eigenvalue $\lambda=1$. Thus all eigenvalues
of $I+JC$ are real and greater than or equal to one. This implies
that matrix $I+JC$ is always invertible, which further implies that $I + CJ = (I+JC)^\top$ is invertible as well.
\end{proof}

\section{Proof of LQT backward conditional HJB equation}\label{sec:Proof_backward_LQT_append}

In this appendix, we prove Theorem \ref{thm:Backward_HJB_LQT}. For LQT, the conditional value function $V(x,s;z,\tau)$ is of the
form \eqref{eq:LQT_V_form1} and \eqref{eq:LQT_V_form2}. At the maximum $\lambda_{*}$ in \eqref{eq:LQT_V_form1}, we have
\begin{align}
C\left(s,\tau\right)\lambda_{*} & =-\left(z-A\left(s,\tau\right)x-b\left(s,\tau\right)\right),\label{eq:lambda_opt}
\end{align}
which implies that $\lambda_{*}$ depends on $x$, $z$, $C\left(s,\tau\right)$, $A\left(s,\tau\right)$, and $b\left(s,\tau\right)$.
We can then write the conditional value function as
\begin{alignat}{1}
& V(x,s;z,\tau)\nonumber \\
& =\zeta+\frac{1}{2}x^{\top}J\left(s,\tau\right)x-x^{\top}\eta\left(s,\tau\right)\nonumber \\
& -\frac{1}{2}\lambda_{*}^{\top}C\left(s,\tau\right)\lambda_{*}-\lambda_{*}^{\top}\left(z-A\left(s,\tau\right)x-b\left(s,\tau\right)\right)\label{eq:V_cond_back_LQT_app}
\end{alignat}
where $\lambda_{*}$ satisfies (\ref{eq:lambda_opt}), and we have dropped the explicit dependence of $\zeta$ on $s$ and $\tau$ for notational simplicity.

In the following, we prove Theorem \ref{thm:Backward_HJB_LQT} by calculating the left-hand side and right-hand
side of \eqref{eq:backward_conditional_HJB} for (\ref{eq:V_cond_back_LQT_app}) and by identifying the terms.

\subsection{Left-hand side of \eqref{eq:backward_conditional_HJB}}

Substituting (\ref{eq:V_cond_back_LQT_app}) into the left-hand side
of \eqref{eq:backward_conditional_HJB} yields
\begin{align}
-\frac{\partial V(x,s;z,\tau)}{\partial s} & =-\frac{\partial \zeta}{\partial s}  -\frac{1}{2}x^{\top}\frac{\partial J\left(s,\tau\right)}{\partial s}x+x^{\top}\frac{\partial\eta\left(s,\tau\right)}{\partial s}\nonumber \\
& +\frac{\partial\lambda_{*}^{\top}}{\partial s}C\left(s,\tau\right)\lambda_{*}+\frac{1}{2}\lambda_{*}^{\top}\frac{\partial C\left(s,\tau\right)}{\partial s}\lambda_{*}\nonumber \\
& +\frac{\partial\lambda_{*}^{\top}}{\partial s}\left(z-A\left(s,\tau\right)x-b\left(s,\tau\right)\right)\nonumber \\
& +\lambda_{*}^{\top}\left(-\frac{\partial A\left(s,\tau\right)}{\partial s}x-\frac{\partial b\left(s,\tau\right)}{\partial s}\right).
\end{align}
Substituting (\ref{eq:lambda_opt}) into the previous equation, we
obtain
\begin{align}
-\frac{\partial V(x,s;z,\tau)}{\partial s} & =-\frac{\partial \zeta}{\partial s} -\frac{1}{2}x^{\top}\frac{\partial J\left(s,\tau\right)}{\partial s}x+x^{\top}\frac{\partial\eta\left(s,\tau\right)}{\partial s}\nonumber \\
& +\frac{1}{2}\lambda_{*}^{\top}\frac{\partial C\left(s,\tau\right)}{\partial s}\lambda_{*}\nonumber \\
& +\lambda_{*}^{\top}\left(-\frac{\partial A\left(s,\tau\right)}{\partial s}x-\frac{\partial b\left(s,\tau\right)}{\partial s}\right).\label{eq:LHS_backward_derivative}
\end{align}

\subsection{Right-hand side of \eqref{eq:backward_conditional_HJB}}

We first calculate
\begin{align}
\frac{\partial V(x,s;z,\tau)}{\partial x} & =J\left(s,\tau\right)x-\eta\left(s,\tau\right)-\frac{\partial\lambda_{*}^{\top}}{\partial x}C\left(s,\tau\right)\lambda_{*}\nonumber \\
& +A^{\top}\left(s,\tau\right)\lambda_{*}\nonumber \\
& -\frac{\partial\lambda_{*}^{\top}}{\partial x}\left(z-A\left(s,\tau\right)x-b\left(s,\tau\right)\right),
\end{align}
where we have used (\ref{eq:V_cond_back_LQT_app}). 

Substituting (\ref{eq:lambda_opt}) into the previous equation, we
obtain
\begin{align}
\frac{\partial V\left(x,s;z,\tau\right)}{\partial x} & =J\left(s,\tau\right)x-\eta\left(s,\tau\right)+A^{\top}\left(s,\tau\right)\lambda_{*}.
\end{align}
Then, the right-hand side of \eqref{eq:backward_conditional_HJB} can be written as
\begin{align}
& \min_{u}\left[\frac{1}{2}(r-Hx)^{\top}X(r-Hx)+\frac{1}{2}u^{\top}Uu\right.\nonumber \\
& \left.+\left(x^{\top}J\left(s,\tau\right)-\eta^{\top}\left(s,\tau\right)+\lambda_{*}^{\top}A\left(s,\tau\right)\right)\left(Fx+Lu+c\right)\right]\label{eq:RHS_backward_1}
\end{align}
where we have removed the dependence of $r$, $H$, $U$, $X$, $F$,
$L$, and $c$ on $s$ for notational simplicity.

Calculating the derivative of (\ref{eq:RHS_backward_1}) with respect to $u$
and setting it to zero, the minimum is achieved
for
\begin{align}
u & =-U^{-1}L^{\top}J^{\top}\left(s,\tau\right)x\nonumber \\
& \quad+U^{-1}L^{\top}\eta\left(s,\tau\right)-U^{-1}L^{\top}A^{\top}\left(s,\tau\right)\lambda_{*}.\label{eq:u_min_backward}
\end{align}
To obtain a solution that ensures that $J\left(s,\tau\right)$ is
a symmetric matrix, we use the identity
\begin{align}
x^{\top}Ax & =\frac{1}{2}x^{\top}\left(A+A^{\top}\right)x\label{eq:symmetric_matrix_identity_append}
\end{align}
which holds for any square matrix $A$. We substitute (\ref{eq:u_min_backward})
into (\ref{eq:RHS_backward_1}), and then use (\ref{eq:symmetric_matrix_identity_append}).
After a lengthy derivation, we obtain a formula of the form (\ref{eq:LHS_backward_derivative}).
By comparing this formula with (\ref{eq:LHS_backward_derivative}),
we identify the derivatives in Theorem \ref{thm:Backward_HJB_LQT} finishing its proof. This appendix has also proved that the conditional value function for the LQT problem is of the form \eqref{eq:LQT_V_form1}.

\section{Proof of LQT forward conditional HJB equation} \label{sec:Proof_forward_LQT_append}

In this appendix, we prove Theorem \ref{thm:Forward_HJB_LQT}. We proceed as in Appendix \ref{sec:Proof_backward_LQT_append}, but using the forward conditional HJB
equation \eqref{eq:forward_conditional_HJB} instead of the backward version \eqref{eq:backward_conditional_HJB}.

\subsection{Left-hand side of \eqref{eq:forward_conditional_HJB}}

Calculating the derivative of (\ref{eq:V_cond_back_LQT_app}) with respect to
$\tau$ and using (\ref{eq:lambda_opt}), we obtain
\begin{align}
\frac{\partial V(x,s;z,\tau)}{\partial\tau} & =\frac{\partial \zeta}{\partial \tau}+\frac{1}{2}x^{\top}\frac{\partial J\left(s,\tau\right)}{\partial\tau}x-x^{\top}\frac{\partial\eta\left(s,\tau\right)}{\partial\tau}\nonumber \\
& -\frac{1}{2}\lambda_{*}^{\top}\frac{\partial C\left(s,\tau\right)}{\partial\tau}\lambda_{*}\nonumber \\
& +\lambda_{*}^{\top}\left(\frac{\partial A\left(s,\tau\right)}{\partial\tau}x+\frac{\partial b\left(s,\tau\right)}{\partial\tau}\right).\label{eq:LHS_forward_derivative}
\end{align}

\subsection{Right-hand side of \eqref{eq:forward_conditional_HJB}}

We first calculate
\begin{align}
\frac{\partial V(x,s;z,\tau)}{\partial z} & =-\frac{\partial\lambda_{*}^{\top}}{\partial z}C\left(s,\tau\right)\lambda_{*}-\lambda_{*}\nonumber \\
& \quad-\frac{\partial\lambda_{*}^{\top}}{\partial z}\left(z-A\left(s,\tau\right)x-b\left(s,\tau\right)\right),
\end{align}
where we have used (\ref{eq:V_cond_back_LQT_app}). Substituting (\ref{eq:lambda_opt})
into the previous expression yields
\begin{align}
\frac{\partial V(x,s;z,\tau)}{\partial z} & =-\lambda_{*}.
\end{align}
Then, the right-hand side of \eqref{eq:forward_conditional_HJB} becomes
\begin{align}
& \min_{u}\left[\frac{1}{2}(r-Hz)^{\top}X(r-Hz)\right.\nonumber \\
& \quad\left.+\frac{1}{2}u^{\top}Uu+\lambda_{*}^{\top}\left(Fz+Lu+c\right)\right], \label{eq:RHS_forward_1}
\end{align}
where we have removed the dependence of $r$, $H$, $U$, $X$, $F$,
$L$, and $c$ on $\tau$ for notational simplicity.

Calculating the derivative of (\ref{eq:RHS_forward_1}) with respect to $u$
and setting it to zero, we obtain
\begin{align}
u & =-U^{-1}L^{\top}\lambda_{*}.\label{eq:u_min_forward}
\end{align}
Substituting (\ref{eq:u_min_forward}) into (\ref{eq:RHS_forward_1}),
we obtain
\begin{align}
& \frac{1}{2}r^{\top}Xr-z^{\top}H^{\top}Xr+\frac{1}{2}z^{\top}H^{\top}XHz\nonumber \\
& +\frac{1}{2}\lambda_{*}^{\top}LU^{-1}L^{\top}\lambda_{*}\nonumber \\
& +\lambda_{*}^{\top}\left(Fz-LU^{-1}L^{\top}\lambda_{*}+c\right). \label{eq:RHS_forward2_append}
\end{align}
This expression is in terms of $z$. To compare with (\ref{eq:LHS_forward_derivative})
and identify the terms, we need to have the expression in terms of
$x$. From (\ref{eq:lambda_opt}), we can write
\begin{align}
z & =A\left(s,\tau\right)x+b\left(s,\tau\right)-C\left(s,\tau\right)\lambda_{*}.
\end{align}
We substitute this expression into (\ref{eq:RHS_forward2_append})
and then use (\ref{eq:symmetric_matrix_identity_append}). After a
lengthy derivation, we obtain an expression of the form (\ref{eq:LHS_forward_derivative}),
which can be used to identify the derivatives in Theorem \ref{thm:Forward_HJB_LQT}, finishing
its proof. This appendix has also proved that the conditional value function for LQT is of the form \eqref{eq:LQT_V_form1}.

\section{Interpolation for non-linear control} \label{sec:Interpolation_append}
This appendix explains the quadratic interpolation method we use in the one-dimensional nonlinear control problem in Section \ref{subsec:nonlinear_control_experiment}. To compute the combination rule in Theorem 2 the direct approach in a nonlinear setting is to use a grid and select the minimum according to the grid. This approach can lead to numerical inaccuracies so we provide the following interpolation method to improve the accuracy. 

Let us consider the minimisation of a general function
\begin{equation}
\begin{split}
g(x) + f(x),
\end{split}
\end{equation}
where we only know $g_{i}=g(x_{i})$ and $f_{i}=f(x_{i})$ at discrete values $x_{i}$, such that $\Delta=x_{i}-x_{i-1}$ for all $i$. We can now find the index $i$ that minimises the sum. We assume that the true minimum lies in $[x_{i-1},x_{i+1}]$. To find a better approximation to the true minimum, we form the interpolants
\begin{equation}
\begin{split}
L_g(x_i + s) &= a_g s^2 + b_g s + c_g, \\
L_f(x_i + s) &= a_f s^2 + b_f s + c_f.
\end{split}
\end{equation}
To obtain the parameters of the interpolants, we denote $g_{i-1}=g(x_{i-1}),\,g_{i}=g(x_{i}),\,g_{i+1}=g(x_{i+1})$. Then, the following equations hold 
\begin{equation}
\begin{split}
g_{i-1} &= a_g \Delta^2 - b_g \Delta + c_g, \\
g_i &= c_g, \\
g_{i+1} &= a_g \Delta^2 + b_g \Delta + c_g.
\end{split}
\end{equation}

Solving these equations, the parameters of the interpolant of $g(\cdot)$ are
\begin{equation}
\begin{split}
c_g &= g_i  \\
b_g &= (g_{i+1} - g_{i-1}) / (2 \Delta) \\
a_g &= (g_{i+1} + g_{i-1} - 2 g_i) / (2 \Delta^2).
\end{split}
\end{equation}
The analogous result can be obtained for the interpolant of $f(\cdot)$. As the middle point is the minimum of the three, we have that $a_{g}+a_{f}\geq 0$. If $a_{g}+a_{f}>0$, the minimum of $L_g(x_i + s)+L_f(x_i + s)$ can be directly obtained as
\begin{equation}
g_i+f_i-\frac{(b_{g}+b_{f}){}^{2}}{4(a_{g}+a_{f})}.
\end{equation}
We can see that the original minimum $g_i+f_i$ is corrected downwards by the interpolation. If $a_{g}+a_{f}=0$, then the minimum is $g_i+f_i$. 

\bibliographystyle{IEEEtran}
\bibliography{IEEEabrv,hjb}

\begin{IEEEbiography}[{\includegraphics[width=1in,height=1.25in,clip,keepaspectratio]{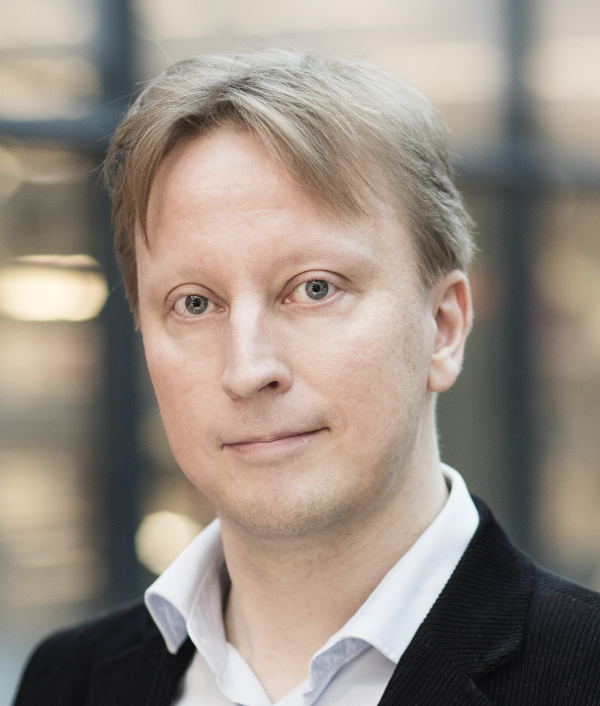}}]{Simo S\"arkk\"a} received his MSc.\ degree (with distinction) in engineering physics and mathematics, and DSc.\ degree (with distinction) in electrical and communications engineering from Helsinki University of Technology, Espoo, Finland, in 2000 and 2006, respectively. Currently, he is a Professor with Aalto University. His research interests are in multi-sensor data processing and control systems. He has authored or coauthored over 200 peer-reviewed scientific articles and two books. He is a Senior Member of IEEE.  \end{IEEEbiography}

\begin{IEEEbiography}[{\includegraphics[width=1in,height=1.25in,clip,keepaspectratio]{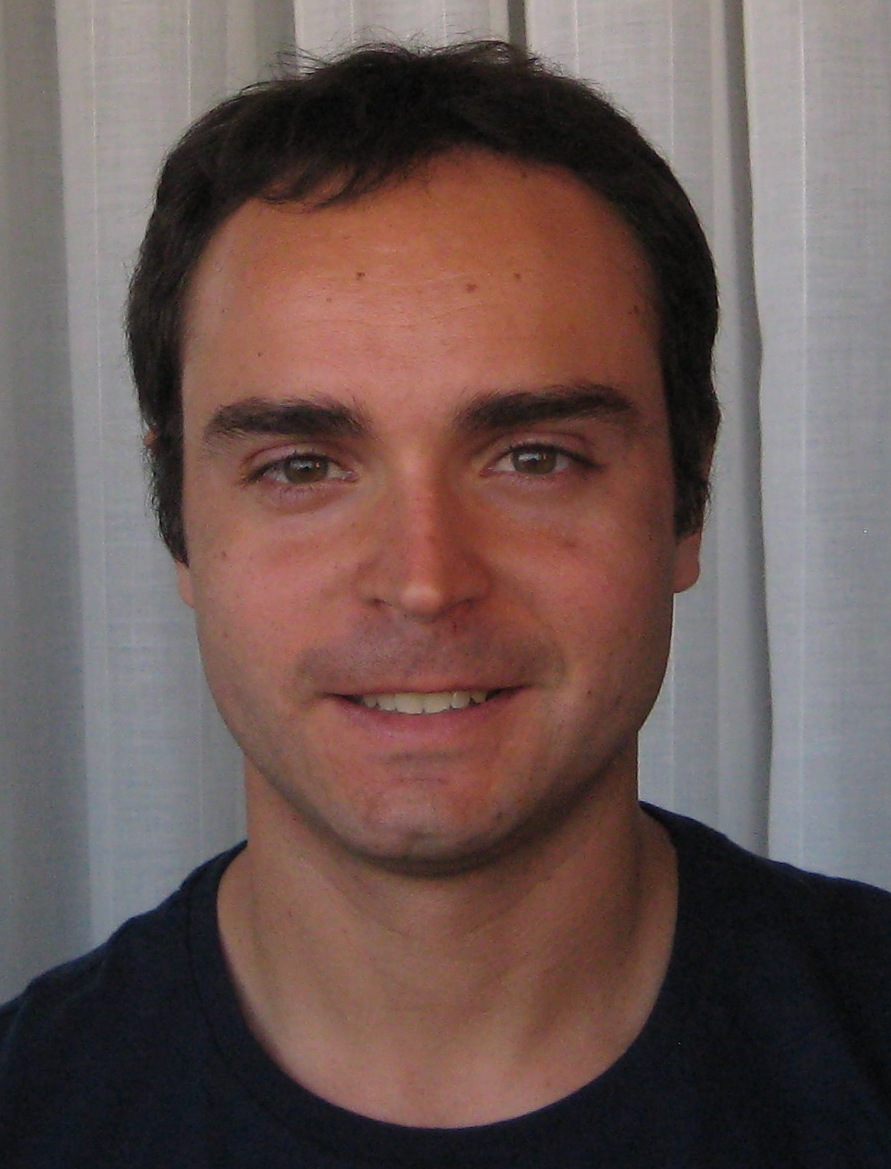}}]{\'Angel F. Garc\'ia-Fern\'andez} received the telecommunication engineering degree (with honours) and the Ph.D. degree from Universidad Polit\'ecnica de Madrid, Madrid, Spain, in 2007 and 2011, respectively. \\ 
He is currently a Senior Lecturer in the Department of Electrical Engineering and Electronics at the University of Liverpool, Liverpool, UK. His main research activities and interests are in the area of Bayesian  estimation, with emphasis on dynamic systems and multiple target tracking. He was the recipient of paper awards at the International Conference on Information Fusion in 2017, 2019 and 2021. \end{IEEEbiography}

\vfill

\end{document}